\documentstyle[12pt]{article}

\newlength{\abstractwidth}
\renewcommand{\thefootnote}{\fnsymbol{footnote}}
\renewcommand{\thanks}[1]{\footnote{#1}} 
\newcommand{\starttext}{ \setcounter{footnote}{0}
\renewcommand{\thefootnote}{\arabic{footnote}}}

\newcommand{\be}{\begin{equation}}
\newcommand{\bea}{\begin{eqnarray}}
\newcommand{\eea}{\end{eqnarray}} \newcommand{\ee}{\end{equation}}

\def\ba{\begin{eqnarray}}
\def\ea{\end{eqnarray}}

\def\o{\omega}
\def\Re{{\rm Re}}

\def\tr{{\rm tr}}
\def\det{{\rm det}}

\def\log{\,{\rm log}\,}

\def\o{\omega}

\def\e{\varepsilon}

\def\o{\omega}
\def\f{\phi}

\def\na{\nabla}

\def\R{{\bf R}}
\def\C{{\bf C}}

\def\p{\prod}

\def\F{{\cal F}}

\def\na{{\nabla}}

\def\[{{\bf [}}
\def\]{{\bf ]}}

\def\p{\partial}

\def\U{\underline{u}}

\begin{document}
\starttext \baselineskip=15pt \setcounter{footnote}{0}
\newtheorem{theorem}{Theorem}
\newtheorem{lemma}{Lemma}
\newtheorem{definition}{Definition}
\newtheorem{proposition}{Proposition}
\newtheorem{corollary}{Corollary}

\begin{center}
{\Large \bf  FULLY NON-LINEAR PARABOLIC EQUATIONS ON COMPACT HERMITIAN MANIFOLDS
\footnote{The first author was supported in part by the National Science Foundation under NSF Grant DMS-12-66033. The second author was supported by the CFM foundation and ATUPS travel grant. }}
\bigskip\bigskip

\centerline{Duong H. Phong and Dat  T. T\^o}

\end{center}

\medskip

\begin{abstract}

{\small 
A notion of parabolic $C$-subsolutions is introduced for parabolic equations, extending the theory of $C$-subsolutions recently developed by B. Guan and more specifically G. Sz\'ekelyhidi for elliptic equations. The resulting parabolic theory provides a convenient unified approach for the study of many geometric flows.}

\end{abstract}

\section{Introduction}
\setcounter{equation}{0}

Subsolutions play an important role in the theory of partial differential equations. Their existence can be viewed as an indication of the absence of any global obstruction.  Perhaps more importantly, it
can imply crucial a priori estimates, as for example in the Dirichlet problem for the complex Monge-Amp\`ere equation \cite{Sp, Gb1}. However, for compact manifolds without boundary, it is necessary to extend the notion of subsolution, since the standard notion may be excluded by either the maximum principle or cohomological constraints. Very recently, more flexible and compelling notions of subsolutions have been proposed by Guan \cite{Gb2} and Sz\'ekelyhidi \cite{Sze}. In particular, they show that their notions, called C-subsolution in \cite{Sze}, do imply the existence of solutions and estimates for a wide variety of fully non-linear elliptic equations on Hermitian manifolds. It is natural to consider also the parabolic case. This was done by Guan, Shi, and Sui in \cite{GSS} for the usual notion of subsolution
and for the Dirichlet problem. 
We now carry this out for the more general notion of $C$-subsolution on compact Hermitian manifolds, adapting the methods of \cite{Gb2} and especially \cite{Sze}. 
As we shall see, the resulting parabolic theory provides a convenient unified approach to the many parabolic equations which have been studied in the literature.

\medskip
Let $(X,\alpha)$ be a compact Hermitian manifold of dimension $n$, $\alpha=i\,\alpha_{\bar{k} j}dz^j\wedge d\bar{z}^k>0$,
and $\chi(z)$ be a real $(1,1)$- form, $$\chi=i\,\chi_{\bar kj}(z) dz^j\wedge d\bar z^k .$$
If $u\in C^2(X)$, let  $A[u]$ be the matrix with entries $A[u]^k{}_j=\alpha^{k\bar{m}}(\chi_{\bar m j}+\p_j\p_{\bar m}u)$. We consider the fully nonlinear parabolic equation,
\begin{equation} \label{parabolic eq}
 \p_tu=F(A[u])- \psi(z),
\end{equation} 
where $F(A)$  is a smooth symmetric function $F(A)=f(\lambda[u])$ of 
the eigenvalues $\lambda_j[u]$, $1\leq j\leq n$ of $A[u]$, 
defined  on a open symmetric, convex cone $\Gamma \subset \R^n$ with vertex at the origin and containing the positive orthant $\Gamma_n$. We shall assume throughout the paper that $f$
satisfies the following conditions:

\smallskip

(1) $f_i>0$ for all $i$, and $f$ is concave. \label{condition 1}

(2) $f(\lambda)\rightarrow -\infty$ as $\lambda\rightarrow \partial \Gamma$ 

(3) For any $\sigma<\sup _\Gamma f$ and $\lambda \in \Gamma$, we have $\lim_{t\rightarrow \infty}f(t\lambda)>\sigma$. \label{condition 3}

 \medskip

We shall say that a $C^2$ function $u$ on $X$ is admissible if the vector of eigenvalues of the corresponding matrix $A$ is in $\Gamma$ for any $z\in X$. 
Fix $T\in (0,\infty]$.  To alleviate the terminology, we shall also designate by the same adjective functions in $C^{2,1}(X\times [0,T))$ which are admissible for each fixed $t\in [0,T)$. The following notion of subsolution is an adaptation to the parabolic case of Sz\'ekelyhidi's \cite{Sze} notion in the elliptic case:
 
\begin{definition}
\label{subsolution}
 An admissible function $\underline{u} \in C^{2,1}(X\times [0,T))$ is said to be a (parabolic) $C$-subsolution of ({\rm \ref{parabolic eq}}), if there exist constants $\delta, K>0$, so that for any $(z,t)\in X\times [0,T)$, the condition
\bea
\label{sub0}
f(\lambda[\U(z,t)]+\mu)-\p_t\U+\tau=\psi(z),
\quad
\mu+\delta I\in \Gamma_n,
\quad
\tau>-\delta
\eea
implies that $|\mu|+|\tau|<K$. Here $I$ denotes the vector $(1,\cdots,1)$ of eigenvalues of the identity matrix.
\end{definition}

\medskip

We shall see below (\S  4.1) that this notion is more general than the classical notion defined by
$f(\lambda([\underline{u}]))-\p_t\U(z,t)>\psi(z,t)$ and studied by Guan-Shi-Sui \cite{GSS}. A $C$-subsolution in the sense of Sz\'ekelyhidi of the equation $F(A[u])-\psi=0$ can be viewed as a parabolic $C$-subsolution of the equation (\ref{parabolic eq}) which is time-independent. But more generally, to solve the equation $F(A[u])-\psi=0$ by say the method of continuity, we must choose a time-dependent deformation of this equation, and we would need then a $C$-subsolution for each time. The heat equation (\ref{parabolic eq}) and the above notion of parabolic subsolution can be viewed as a canonical choice of deformation.

\medskip

To discuss our results, we need a finer classification of non-linear partial differential operators due to Trudinger \cite{Tr95}. Let $\Gamma_\infty$ be the projection of $\Gamma_n$ onto ${\bf R}^{n-1}$,
\bea
\Gamma_\infty
=
\{\lambda'=(\lambda_1,\cdots,\lambda_{n-1}); \ \lambda=(\lambda_1,\cdots,\lambda_n)\in \Gamma \ {\rm for \ some}\ \lambda_n\}
\eea
and define the function $f_\infty$ on $\Gamma_\infty$ by
\bea
\label{finfinity}
f_\infty(\lambda')={\rm lim}_{\lambda_n\to \infty}f(\lambda',\lambda_n).
\eea
It is shown in \cite{Tr95} that, as a consequence of the concavity of $f$, the limit is either finite for all $\lambda'\in \Gamma_\infty$ or infinite for all $\lambda'\in\Gamma_\infty$. We shall refer to the first case as the {\it bounded case}, and to the second case as the {\it unbounded case}. For example, Monge-Amp\`ere flows belong to the unbounded case, while the $J$-flow and Hessian quotient flows belong to the bounded case. In the unbounded case, any admissible function, and in particular $0$ if $\lambda[\chi]\in\Gamma$, is a $C$-subsolution in both the elliptic and parabolic cases. 
We have then:

\begin{theorem}
\label{theorem B}
Consider the flow {\rm (\ref{parabolic eq})}, and assume that $f$ is in the unbounded case. Then for any 
admissible initial data $u_0$, the flow admits a smooth solution $u(z,t)$ 
on $[0,\infty)$, and its normalization $\tilde u$ defined by
\bea
\label{normalization}
\tilde{u}:= u-{1\over V}\int_X u\, \alpha^n,
\qquad
V=\int_X\alpha^n,
\eea
converges in $C^\infty$ to a function $\tilde u_\infty$ satisfying the following equation for some constant $c$,
\bea
\label{equation}
F(A[\tilde u_\infty])=\psi(z)+c.
\eea
\end{theorem}

\medskip
The situation is more complicated when $f$ belongs to the bounded case:
\medskip

\begin{theorem}
\label{theorem C}
Consider the flow {\rm (\ref{parabolic eq})}, and assume  that it admits a subsolution $\U$ on $X\times [0,\infty)$, but that $f$ is in the bounded case. Then 
for any admissible data $u_0$, the equation admits a smooth solution $u(z,t)$
on $(0,\infty)$. Let $\tilde u$ be the normalization
of the solution $u$, defined as before by (\ref{normalization}). Assume that
either one of the following two conditions holds.

\smallskip

{\rm (a)} The initial data and the subsolution satisfy 
\bea
\p_t\U
\geq {\rm sup}_X (F(A[u_0])-\psi);
\eea

{\rm (b)} or there exists a function $h(t)$ with $h'(t)\leq 0$ so that
\bea
{\rm sup}_X (u(t)-h(t)-\underline{u}(t))\geq 0
\eea
and the Harnack inequality
\bea
\label{harnack_ineq}
{\rm sup}_X(u(t)-h(t))\leq -C_1{\rm inf}_X (u(t)-h(t))+C_2
\eea
holds for some constants $C_1,C_2>0$ independent of time.

\smallskip
\noindent
Then $\tilde u$ converges in $C^\infty$ to a function $\tilde u_\infty$ satisfying (\ref{equation}) for some constant $c$.
\end{theorem}

The essence of the above theorems resides in the a priori estimates which are established in \S 2. The $C^1$ and $C^2$ estimates can be adapted from the corresponding estimates for $C$-subsolutions in the elliptic case, but the $C^0$ estimate turns out to be more subtle. 
Following Blocki \cite{B} and Sz\'ekelyhidi \cite{Sze}, we obtain $C^0$ estimates from the Alexandrov-Bakelman-Pucci (ABP) inequality, using this time a parabolic version of ABP due to K. Tso \cite{Ts}. 
However, it turns out that the existence of a $C$-subsolution gives only partial information on the oscillation of $u$, and what can actually be estimated has to be formulated with some care, leading to the distinction between the cases of $f$ bounded and unbounded, as well as Theorem \ref{theorem C}.

\smallskip
The conditions (a) and especially (b) in Theorem 2 may seem impractical at first sight since they involve the initial data as well as the long-time behavior of the solution. Nevertheless, as we shall discuss in greater detail in section \S 4, Theorems 1 and 2 can be successfully applied to a wide range of parabolic flows on Hermitian manifolds previously studied in the literature, including the K\"ahler-Ricci flow, the Chern-Ricci flow, the $J$-flow, the Hessian flows, the quotient Hessian flows, and mixed Hessian flows. We illustrate this by deriving in \S 4 as a corollary of Theorem 2  a convergence theorem for a mixed Hessian flow, which 
seems new to the best of our knowledge. It answers a question raised for general $1\leq \ell <k\leq n$ by Fang-Lai-Ma \cite{FLM} (see also Sun \cite{Su1, Su2, Su3, Su5}), and extends the solution obtained for $k=n$ by Collins-Sz\'ekelyhidi \cite{CSze} and subsequently also by Sun \cite{Su5, Su6}:

\begin{theorem}
\label{Cor9}
Assume that $(X,\alpha)$ is a compact K\"ahler $n$-manifold, and fix $1\leq \ell<k\leq n$. Fix a closed $(1,1)$-form $\chi$ which is $k$-positive and non-negative constants $c_j$, and assume that there exists a form $\chi'=\chi+i\p\bar\p \underline{u}$ which is a closed $k$-positive form and satisfies
\begin{equation}\label{positivity_condition_3}
kc(\chi')^{k-1}\wedge\alpha^{n-k}-\sum_{j=1}^{\ell}j c_j(\chi')^{j-1}\wedge \alpha^{n-j}>0,
\end{equation}
in the sense of positivity of $(n-1,n-1)$-forms. Here the constant $c$ is given by
\bea
c[\chi^{k}][\alpha^{n-k}]=\sum_{j=1}^\ell c_j[\chi^{j}][\alpha^{n-j}].
\eea
Then the flow
\bea
\label{combination}
\p_tu
=
-{\sum_{j=1}^\ell c_j\sigma_j(\lambda(A[u]))\over \sigma_k(\lambda(A[u]))}+c,
\qquad u(\cdot,0)=0,
\eea
admits a solution for all time which converges smoothly to a function $u_\infty$ as $t\to\infty$. The
form $\o=\chi+i\p\bar\p u_\infty$ is $k$-positive and satisfies the equation
\bea
c\,\omega^k\wedge \alpha^{n-k}=\sum_{j=1}^\ell c_j\, \omega^j\wedge \alpha^{n-j}.
\eea
\end{theorem}

Regarding the condition (a) in Theorem 2, we note that natural geometric flows whose long-time behavior may be very sensitive to the initial data are appearing increasingly frequently in non-K\"ahler geometry. A prime example is the Anomaly flow, studied in \cite{PPZ2, PPZ5, PPZ6, PPZ7, FHP}. Finally, 
Theorem \ref{theorem C} will also be seen to imply as a corollary 
a theorem of Sz\'ekelyhidi (\cite{Sze}, Proposition 26), and the condition for solvability there will be seen to correspond to condition (a) in Theorem \ref{theorem C}. This suggests in particular that some additional conditions for the convergence of the flow cannot be dispensed with altogether.

\section{A Priori Estimates} 
\setcounter{equation}{0}

\subsection{$C^0$ Estimates}

We begin with the $C^0$ estimates implied by the existence of a $C$-subsolution for the parabolic flow (\ref{parabolic eq}). One of the key results of \cite{Sze} was that the existence of a subsolution in the elliptic case implies a uniform bound for the oscillation of the unknown function $u$. In the parabolic case, we have only the following weaker estimate:

\begin{lemma}
\label{C0}
Assume that the equation {\rm (\ref{parabolic eq})} admits a parabolic $C$-solution 
on $X\times [0,T)$ in the sense of Definition \ref{subsolution}, and that there exists a $C^1$ function $h(t)$ with $h'(t)\leq 0$ and
\bea
\label{ABP0}
{\rm sup}_X(u(\cdot,t)-\underline{u}(\cdot,t)-h(t))\geq 0.
\eea
Then there exists a constant $C$ depending only on $\chi,\alpha$, $\delta$, $\|u_0\|_{C^0}$,
and $\|i\p\bar\p \underline{u}\|_{L^\infty}$
so that
\bea
u(\cdot,t)-\underline{u}(\cdot,t)-h(t)\geq -C\ \  {\rm for\ all}
\ (z,t)\in X\times [0,T).
\eea
\end{lemma}

\medskip
\noindent
{\it Proof.} First, note that by Lemma \ref{pt} proven later in \S 3, the function $\p_t u$ is uniformly bounded for all time by a constant depending only on $\psi$ and the initial data $u_0$. Integrating this estimate on $[0,\delta]$ gives a bound for $|u|$ on $X\times[0,\delta]$ depending only on $\psi$, $u_0$ and $\delta$.
Thus we need only consider the range $t\geq\delta$. Next, the fact that $\U$ is a parabolic subsolution and the condition that $h'(t)\leq 0$ imply that $\U+h(t)$ is a parabolic subsolution as well. So it suffices to prove the desired inequality with $h(t)=0$, as long as the constants involved do not depend on $\p_t\U$. Fix now any $T'<T$, and set for each $t$, $v=u-\U$, and 
\bea
L={\rm min}_{X\times[0,T']} v=v(z_0,t_0)
\eea
for some $(z_0,t_0)\in X\times [0,T']$. We shall show that $L$ can be bounded from below by a constant depending only on the initial data $u_0$ and independent of $T'$. We can assume that $t_0>0$, otherwise we are already done. Let $(z_1,\cdots,z_n)$ be local holomorphic coordinates for $X$ centered at $z_0$, $U=\{z;|z|<1\}$, and define the following function on the set ${\cal U}=U\times \{t;-\delta\leq 2(t-t_0)<\delta\}$,
\bea
w=v+{\delta^2\over 4}|z|^2+|t-t_0|^2,
\eea
where $\delta>0$ is the constant appearing in the definition of subsolutions. Clearly $w$ attains its minimum on ${\cal U}$ at $(z_0,t_0)$, and $w\geq {\min}_{\cal U}w+{1\over 4}\delta^2$ on the parabolic boundary of ${\cal U}$. We can thus apply the following parabolic version of the Alexandrov-Bakelman-Pucci inequality, due to K. Tso (\cite{Ts}, Proposition 2.1, with the function $u$ there set to $u=-w+{\rm min}_{\cal U}w+{\delta^2\over 4}$):

\medskip
Let ${\cal U}$ be the subset of ${\bf R}^{2n+1}$ defined above,
and let  $w: {\cal U}\rightarrow {\bf R}$ be a smooth function which attains its minimum at $(0,t_0)$, and $w\geq {\rm min}_{\cal U}w+{1\over 4}\delta^2$ on the parabolic boundary of ${\cal U}$. Define the set
\begin{equation}
\label{contact set}
S:=\left\lbrace (x,t) \in {\cal U}: \begin{array}{l}
  w(x,t)\leq w(z_0,t_0)+{1\over 4}\delta^2,\quad |D_x w(x,t) |<\frac{\delta^2}{8}, \textit{ and }\\
  w(y,s)\geq w(x,t)+D_x w(x,t). (y-x),\, \forall y\in U, s\leq t 
 \end{array}
\right\rbrace.
\end{equation}
 Then there is a constant $C=C(n)>0$ so that
 $$C\delta^{4n+2}\leq \int_{S} (-w_t)\, \det(w_{ij})dxdt. $$
      
\medskip
Returning to the proof of Lemma \ref{C0}, we claim that, on the set $S$, we have
\bea
\label{ABP1}
|w_t|+{\rm det}\,(D^2_{jk}w)\leq C
\eea
for some constant depending only on $\delta$, and $\|i\p\bar\p \U\|_{L^\infty}$. Indeed, let
\bea
\mu=\lambda[u]-\lambda[\U],
\qquad
\tau=-\p_tu+\p_t\U.
\eea
Along $S$, we have $D^2_{ij}w\geq 0$ and $\p_tw\leq 0$. In terms of $\mu$ and $\tau$, this means that $\mu+\delta {\bf I}\in \Gamma_n$ and $0\leq -\p_tw=\tau-2(t-t_0)\leq \tau+\delta$. The fact that $u$ is a solution of the equation (\ref{parabolic eq}) can be expressed as
\bea
f(\lambda[\U]+\mu)-\p_t\U+\tau=\psi(z).
\eea
Thus the condition that $\U$ is a parabolic subsolution implies that $|\mu|$ and $|\tau|$ are bounded uniformly in $(z,t)$. Since along $S$, we have
${\rm det}(D^2_{ij}w)\leq 2^n( {\rm det}(D^2_{\bar kj}w))^2$, it follows that both $|w_t|$ and ${\rm det}(D^2_{ij}w)$ are bounded uniformly, as was to be shown.

\smallskip
Next, by the definition of the points $(x,t)$ on $S$, we have $w(x,t)\leq L+{\delta^2\over 4}$. Since we can assume that $|L|>\delta^2$, it follows that $w<0$ and $|w|\geq {|L|\over 2}$ on $S$. Thus we can write, in view of (\ref{ABP1}), for any $p>0$,
\bea
\label{ABP2}
C_n\delta^{4n+2}
\leq
C\int_S dxdt
\leq
\left({|L|\over 2}\right)^{-p}\int_S |w(x,t)|^pdxdt
\leq
\left({|L|\over 2}\right)^{-p}\int_{\cal U} |w(x,t)|^pdxdt.
\eea
Next write
\bea
\label{ABP3}
|w|
&=&-w=-v-{\delta^2\over 4} |z|^2-(t-t_0)^2
\leq -v\nonumber\\
&\leq& -v+{\rm sup}_Xv
\eea
since ${\rm sup}_Xv\geq 0$ by the assumption (\ref{ABP0}). Since $\lambda[u]\in \Gamma$ and the cone $\Gamma$ is convex, it follows that $\Delta u\geq -C$ and hence
\bea
\Delta (v-{\rm sup}_X\,v)=\Delta u-\Delta\U\geq -A
\eea
for some constant $A$ depending only on $\chi$, $\alpha$, and $\|i\p\bar\p\U\|_{L^\infty}$. The Harnack inequality applied to the function $v-{\rm sup}_Xv$, in the version provided by Proposition 10, \cite{Sze}, implies that 
\bea
\|v-{\rm sup}_X v\|_{L^p(X)}
\leq C
\eea
for $C$ depending only on $(X,\alpha),A$, and $p$. Substituting these bounds into (\ref{ABP2}) gives
\bea
C\delta^{4n+2}
\leq \left({|L|\over 2}\right)^{-p}\int_{|t|<{1\over 2}\delta}\|{\rm sup}_Xv-v\|_{L^p(X)}^p dt
\leq C'\delta \left({|L|\over 2}\right)^{-p}
\eea
from which the desired bound for $L$ follows. Q.E.D.

\subsection{$C^2$ Estimates}

In this section we prove an estimate for the complex Hessian of $u$ in terms of the gradient. The original strategy goes back to the work of Chou-Wang \cite{CW}, with adaptation to complex Hessian equations by Hou-Ma-Wu \cite{HMW}, and to fully non-linear elliptic equations admitting a $C$-subsolution by Guan \cite{Gb2} and Sz\'ekelyhidi \cite{Sze}. Other adaptations to $C^2$ estimates can be found in \cite{STW15}, \cite{PPZ1}, \cite{PPZ3}, \cite{Zhe16}. We follow closely \cite{Sze}.

\begin{lemma}
\label{C2}
Assume that the flow $({\rm \ref{parabolic eq}})$ admits a $C$-subsolution on $X\times [0,T)$. Then we have the following estimate
\bea
|i\p\bar\p u |\leq \tilde{C} (1+ \sup_{ X\times [0,T)}| \nabla u|^2_\alpha) 
\eea
where $\tilde{C}$ depends only on $\| \alpha\|_{C^2},\|\psi \|_{C^{2}},\|\chi\|_{C^2},
\|\tilde u-\underline{\tilde u}\|_{L^\infty},\|\nabla\underline{u}\|_{L^\infty},
\|i\p\bar\p \underline{u}\|_{L^\infty}, \|\p_t u\|_{L^\infty}$,  $\|\p_t (u-\U) \|_{L^\infty}$, and the dimension $n$.
\end{lemma}

\bigskip
\noindent
{\it Proof.} Let ${\cal L}=-\p_t+F^{k\bar k}\na_k\na_{\bar k}$. Denote $g=\chi+i\p\bar \p u $, then $A[u]^{k}{}_j=\alpha^{k \bar p} g_{\bar p j}$. We would like to apply the maximum principle to the function
\bea
G=\log \lambda_1+\phi(|\na u|^2)+\varphi(\tilde v)
\eea
where $v=u-\U$, $\tilde v$ is the normalization of $v$, $\lambda_1:X\to {\bf R}$ is the largest eigenvalue of the matrix $A[u]$ at each point, and the functions $\phi$ and $\varphi$ will be specified below. Since the eigenvalues of $A[u]$ may not be distinct, we perturb $A[u]$ following the technique of \cite{Sze}, Proposition 13. Thus assume that $G$ attains its maximum on $X\times[0,T']$ at some $(z_0,t_0)$, with $t_0>0$. We choose local complex coordinates, so that $z_0$ corresponds to $0$, and $A[u]$ is diagonal at $0$ with eigenvalues $\lambda_1\geq \cdots\geq \lambda_n$. Let $B=(B^i{}_j)$ be a diagonal matrix with $0=B^1{}_1<B^2{}_2<\cdots<B^n{}_n$ and small constant entries, and set $\tilde A=A-B$. Then at the origin $\tilde A$ has eigenvalues $\tilde\lambda_1=\lambda_1$, $\tilde \lambda_i=\lambda_i-B^i{}_i<\tilde\lambda_1$ for all $i>1$. 

\smallskip
Since all the eigenvalues of $\tilde A$ are distinct, we can define near $0$ the following smooth function $\tilde G$,
\bea
\tilde G=\log\tilde\lambda_1+
\phi(|\na u|^2)+\varphi(\tilde v)
\eea
where
\bea
\phi(t)=-{1\over 2}\log(1-{t\over 2P}),
\quad
P={\rm sup}_{X\times [0,T']}(|\na u|^2+1)
\eea
and, following \cite{STW15}
\bea
\varphi(t)=D_1 e^{-D_2t}
\eea
for some large constants $D_1,D_2$ to be chosen later. Note that
\bea
{1\over 4P}\leq\phi'\leq{1\over 2P},
\qquad \phi''=2(\phi')^2>0.
\eea

The norm $|\na u|^2$ is taken with respect to the fixed Hermitian metric $\alpha$ on $X$, and we shall compute using covariant derivatives $\na$ with respect to $\alpha$. Since the matrix $B^j{}_m$ is constant in a neighborhood of $0$ and since we are using the Chern unitary connection, we have $\na_{\bar k}B^j{}_m=0$.
Our conventions for the curvature and torsion tensors of a Hermitian metric
$\alpha$ are as follows,
\bea
[\na_\beta,\na_\alpha]V^\gamma=
R_{\alpha\beta}{}^\gamma{}_\delta V^\delta
+
T^\delta{}_{\alpha\beta}\na_\delta V^\gamma.
\eea
We also set
\bea
{\cal F}=\sum_i f_i(\lambda[u]).
\eea
An important observation is that there exists a constant $C_1$, depending only on $\|\psi\|_{L^\infty(X)}$ and $\|\p_tu\|_{L^\infty(X\times [0,T))}$ so that
\bea
{\cal F}\geq C_1.
\eea
Indeed it follows from the properties of the cone $\Gamma$ that $\sum_i f_i(\lambda)\geq C(\sigma)$ for each fixed $\sigma$ and $\lambda\in \Gamma^\sigma$. When $\lambda=\lambda[u]$,
$\sigma$ must lie in the range of $\p_tu+\psi$, which is a compact set
bounded by $\|\p_tu\|_{L^\infty(X\times [0,T))}+\|\psi\|_{L^\infty(X)}$, hence our claim.

\subsubsection{Estimate of ${\cal L}(\log \tilde\lambda_1)$}

Clearly
\bea
{\cal L}\log{\tilde \lambda_1}
={1\over\lambda_1}(F^{k\bar k}\tilde\lambda_{1,\bar kk}-
\p_t\tilde\lambda_1)-F^{k\bar k}{|\tilde\lambda_{1,{\bar k}}|^2\over\lambda_1^2}.
\eea
We work out the term $F^{k\bar k}\tilde\lambda_{1,\bar kk}-
\p_t\tilde\lambda_1$ using the flow.
The usual differentiation rules (\cite{Sp}) readily give
\bea
\label{lambda1k}
\tilde\lambda_{1,\bar k}=\na_{\bar k}g_{\bar 11}
\eea
and
\bea
\tilde\lambda_{1,\bar kk}
=
\na_{ k}\na_{\bar k} g_{\bar 11}
+
\sum_{p>1}{|\na_{\bar k}g_{\bar p1}|^2+|\na_{\bar k} g_{\bar 1p}|^2\over \lambda_1-\tilde\lambda_p}
-
\sum_{p>1}{\na_k B^1{}_p\na_{\bar k}g_{\bar p1}+
\na_kB^p{}_1\na_{\bar k}g_{\bar 1p}\over\lambda_1-\tilde\lambda_p}.
\eea
while it follows from the flow that
\bea
\p_t\tilde\lambda_1=\p_t u_{\bar 11}
=
F^{l\bar k,s\bar r}\na_{\bar 1}g_{\bar kl}\na_1 g_{\bar rs}+F^{k\bar k}\na_1\na_{\bar 1} g_{\bar kk}-\psi_{\bar 11}.
\eea
Thus

\bea
F^{k\bar k}\tilde\lambda_{1,\bar kk}-
\p_t\tilde\lambda_1
&=&
F^{k\bar k}(\na_k\na_{\bar k}g_{\bar 11}-\na_1\na_{\bar 1} g_{\bar kk})+
F^{l\bar k,s\bar r}\na_{\bar 1}g_{\bar kl}\na_1 g_{\bar rs}-\psi_{\bar 11}
\nonumber\\
&&
+
F^{k\bar k}\sum_{p>1} \big\{ {|\na_{\bar k}g_{\bar p1}|^2+|\na_{\bar k} g_{\bar 1p}|^2\over \lambda_1-\tilde\lambda_p}
-
{\na_k B^1{}_p\na_{\bar k}g_{\bar p1}+
\na_kB^p{}_1\na_{\bar k}g_{\bar 1p}\over\lambda_1-\tilde\lambda_p}\big\}
\nonumber
\eea
A simple computation gives
\bea
\na_k\na_{\bar k}g_{\bar 11}-\na_1 \na_{\bar 1}g_{\bar kk}
&=&-2\Re (T_{k1}^p\na_{\bar k} g_{\bar p 1})
+T\star \na\chi+R\star\na\bar\na u+T\star T\star \na\bar\na u
\nonumber\\
&\geq&
-2\Re (T_{k1}^p\na_{\bar k} g_{\bar p1})-C_2(\lambda_1+1),
\eea
where $C_2$ depending only on $\|\alpha\|_{C^2}$ and $\|\chi\|_{C^2}$.
We also have

\bea
&&
\sum_{p>1}\big\{{|\na_{\bar k} g_{\bar p1}|^2+|\na_{\bar k} g_{\bar 1p}|^2\over \lambda_1-\tilde\lambda_p}
-
{\na_k B^1{}_p\na_{\bar k}g_{\bar p1}+
\na_kB^p{}_1\na_{\bar k}g_{\bar 1p}\over\lambda_1-\tilde\lambda_p}\big\}
\\
&&
\geq
{1\over 2}
\sum_{p>1}{|\na_{\bar k}g_{\bar p1}|^2+|\na_{\bar k} g_{\bar 1p}|^2\over \lambda_1-\tilde\lambda_p}-C_3
\geq
{1\over 2(n\lambda_1+1)}
\sum_{p>1}|\na_{\bar k} g_{\bar p1}|^2+|\na_{\bar  k} g_{\bar 1p}|^2-C_3,
\nonumber
\eea
where $C_3$ only depends  on the dimension $n$, and the second inequality is due to the fact that $(\lambda_1-\tilde\lambda_p)^{-1}\geq (n\lambda_1+1)^{-1}$, which follows itself from
the fact that $\sum_i\lambda_i\geq 0$ and $B$ was chosen to be small. Thus
\bea
&&
\na_k\na_{\bar k}g_{\bar 11}-\na_1 \na_{\bar 1}g_{\bar kk}
+
\sum_{p>1}\big\{{|\na_{\bar k}g_{\bar p1}|^2+|\na_{\bar k} g_{\bar 1p}|^2\over \lambda_1-\tilde\lambda_p}
-
{\na_k B^1{}_p\na_{\bar k}g_{\bar p1}+
\na_kB^p{}_1\na_{\bar k}g_{\bar 1p}\over\lambda_1-\tilde\lambda_p}\big\}
\nonumber\\
&&
\geq
-2\Re (T_{k1}^p\na_{\bar k} g_{\bar p1})
+{1\over 2(n\lambda_1+1)}
\sum_{p>1}|\na_{\bar k} g_{\bar p1}|^2+|\na_{\bar k}g_{\bar 1p}|^2-C_2(\lambda_1+1)-C_3
\nonumber\\
&&
\geq
-C_4|\na_{\bar k}g_{\bar 11}|-C_5\lambda_1-C_6
\eea
where we have used the positive terms to absorb all the terms $T_{k1}^p\na_{\bar k} g_{\bar p1}$, except for $T_{k1}^1 \na_{\bar k} g_{\bar 11}$ and  $C_4,C_5,C_6$ only depend on $\|\alpha\|_{C^2}, \|\chi\|_{C^2},n$. Altogether,
\bea
F^{k\bar k}\tilde\lambda_{1,\bar kk}
-
\p_t\tilde\lambda_1
\geq
-C_4F^{k\bar k}|\na_{\bar k}g_{\bar 11}|
+F^{l\bar k,s\bar r}\na_{\bar 1}g_{\bar kl}\na_1 g_{\bar rs}-\psi_{\bar 11}- C_5{\cal F}\lambda_1 - C_6{\cal F}
\eea
and we find
\bea
{\cal L}\log \tilde\lambda_1
&\geq &
-F^{k\bar k}{|\tilde\lambda_{1, \bar k}|^2\over\lambda_1^2}
-{1\over\lambda_1}
F^{l\bar k,s\bar r}\na_{\bar 1}g_{\bar kl}\na_1 g_{\bar rs}
-C_4{1\over\lambda_1}F^{k\bar k}|\na_{\bar k}g_{\bar 11}|
-C_7{\cal F},
\eea
where we have bounded $\psi_{\bar 11}$ by a constant that can be absorbed in $C_6{\cal F}/\lambda_1\leq C_6{\cal F}$, since $\lambda_1\geq 1$ by assumption, and ${\cal F}$ is bounded below by a constant depending on $\|\psi \|_{L^\infty} $ and $\|\p_t u\|_{L^\infty}$. The constant $C_7$ thus only depends on $\|\alpha\|_{C^2},\|\chi\|_{C^2},n$, $\| \p_t u\|_{L^\infty}$ and $\| \psi\|_{C^2}$. In view of (\ref{lambda1k}), this can also be rewritten as
\bea
{\cal L}\log \tilde\lambda_1
&\geq &
-F^{k\bar k}{|\tilde\lambda_{1, \bar k}|^2\over\lambda_1^2}
-{1\over\lambda_1}
F^{l\bar k,s\bar r}\na_{\bar 1}g_{\bar kl}\na_1 g_{\bar rs}
-C_4{1\over\lambda_1}F^{k\bar k}|\tilde\lambda_{1,\bar k}|
-C_7{\cal F}.
\eea

\subsubsection{Estimate for ${\cal L}\phi(|\na u|^2)$}

Next, a direct calculation gives
\bea
{\cal L}\phi(|\na u|^2)
&=&
\phi'(F^{q\bar q}\na_q\na_{\bar q}-\p_t)|\na u|^2
+
\phi''F^{q\bar q}\na_q|\na u|^2\na_{\bar q}|\na u|^2
\nonumber\\
&=&
\phi'\big\{\na^ju(F^{q\bar q}\na_q\na_{\bar q}-\p_t)\na_ju
+
\na^{\bar j}u(F^{q\bar q}\na_q\na_{\bar q}-\p_t)\na_{\bar j}u\big\}
\nonumber\\
&&
+
\phi'
F^{q\bar q}(|\na_q\na u|^2+|\na_q\bar\na u|^2)
+
\phi''F^{q\bar q}\na_q|\na u|^2\na_{\bar q}|\na u|^2.
\eea
In view of the flow, we have
\bea
\na_j\p_tu=F^{k\bar k}\na_j g_{\bar kk}-\psi_j,
\quad
\na_{\bar j}\p_tu=F^{k\bar k}\na_{\bar j} g_{\bar kk}-\psi_{\bar j}.
\eea
It follows that
\bea
(F^{k\bar k}\na_k\na_{\bar k}-\p_t)\na_{\bar j}u&=&
F^{k\bar k}(\na_k\na_{\bar k}u_{\bar j}-\na_{\bar j}g_{\bar kk})+\psi_{\bar j}
\nonumber\\
&=&F^{k\bar k}(-\na_{\bar j}\chi_{\bar kk}+
\overline{T_{kj}^p}\na_{\bar j}\na_ku+R_{\bar jk}{}^{\bar m}{}_{\bar k}\na_{\bar m}u)
+\psi_{\bar j}
\eea
and hence, for small $\varepsilon$, there is a constant $C_8>0$ depending only on $\varepsilon,\|\chi\|_{C^2}, \|\alpha\|_{C^2}$ and $| | \psi\|_{C^2} $ such that
\bea
\phi'\na^{\bar j}u(F^{q\bar q}\na_q\na_{\bar q}-\p_t)\na_{\bar j}u
\geq
-C_8{\cal F}-{\varepsilon\over P}
F^{q\bar q}(|\na_q\na u|^2+|\na_q\bar\na u|^2)
\eea
since we can assume that $\lambda_1>> P={\rm sup}_{X\times [0,T']}(|\na u|^2+1)$
(otherwise the desired estimate $\lambda_1<CP$ already holds), and $(4P)^{-1}<\phi'<(2P)^{-1}$. Similarly we obtain the same estimate for $ \phi' \na^ju(F^{q\bar q}\na_q\na_{\bar q}-\p_t)\na_ju $. Thus by choosing $\varepsilon=1/24$, we have
\bea
{\cal L}\phi(|\na u|^2)
\geq 
-C_8{\cal F}+{1\over 8P}F^{q\bar q}(|\na_q\na u|^2+|\na_q\bar\na u|^2)
+
\phi''F^{q\bar q}\na_q|\na u|^2\na_{\bar q}|\na u|^2.
\eea

\subsubsection{Estimate for ${\cal L}\tilde G$}

The evaluation of the remaining term ${\cal L}\varphi(\tilde v)$ is straightforward,
\bea
{\cal L}\varphi(\tilde v)
=
\varphi'(\tilde v)(F^{k\bar k}\na_k\na_{\bar k}\tilde v-\p_t\tilde v)
+
\varphi''(\tilde v) F^{k\bar k}\na_k\tilde v\na_{\bar k}\tilde v.
\eea
Altogether, we have established the following lower bound for ${\cal L}\tilde G$,
\bea
\label{lower bound}
{\cal L}\tilde G
&\geq& 
-F^{k\bar k}{|\tilde\lambda_{1,\bar k}|^2\over\lambda_1^2}
-{1\over\lambda_1}
F^{l\bar k,s\bar r}\na_{\bar 1}g_{\bar kl}\na_1 g_{\bar rs}
-C_4{1\over\lambda_1}F^{k\bar k}|\lambda_{1,\bar k}|
-C_9{\cal F}\nonumber\\
&&
+{1\over 8P}F^{q\bar q}(|\na_q\na u|^2+|\na_q\bar\na u|^2)
+
\phi''F^{q\bar q}\na_q|\na u|^2\na_{\bar q}|\na u|^2
\nonumber\\
&&
+\varphi'(\tilde v)(F^{k\bar k}\na_k\na_{\bar k}\tilde v-\p_t\tilde v)
+
\varphi''(\tilde v) F^{k\bar k}\na_k\tilde v\na_{\bar k}\tilde v,
\eea
where $C_4$ and $C_9$ only depend on $\|\chi\|_{C^2},\|\alpha\|_{C^2}, \| \psi\|_{C^2}, \|\p_t u \|_{L^\infty}$  and the dimension $n$.

\medskip
\noindent
For a small $\theta>0$ to be chosen hereafter, we deal with two following  cases.

\subsubsection{Case 1: $\theta\lambda_1\leq -\lambda_n$}

In this case, we have $\theta^2\lambda_1^2\leq \lambda_n^2$. Thus we can write
\bea
{1\over 8P}
F^{q\bar q}(|\na_q\na u|^2+|\na_q\bar\na u|^2)
&\geq& {F^{n\bar n}\over 8P}|u_{\bar nn}|^2
=
{F^{n\bar n}\over 8P}|\lambda_n-\chi_{\bar nn}|^2
\geq {{\cal F}\lambda_n^2\over 10n P}-{C_{10}{\cal F}\over P}
\nonumber\\
&\geq&
{\theta^2\over 10 nP}{\cal F}\lambda_1^2-C_{10}{\cal F},
\eea
where $C_{10}$ only depends on $\|\chi\|_{C^2}$. 
Next, it is convenient to combine the first and third terms in the expression for ${\cal L}\tilde G$,
\bea
-F^{k\bar k}{|\tilde\lambda_{1,\bar k}|^2\over\lambda_1^2}
-C_4{1\over \lambda_1}F^{k\bar k}|\tilde \lambda_{1,\bar k}|
\geq -{3\over 2}F^{k\bar k}{|\tilde\lambda_{1,\bar k}|^2\over\lambda_1^2}
-C_{11}{\cal F}.
\eea
where $C_{11}$ only depends on $C_4$.  

\medskip
\noindent
At a maximum point for $\tilde G$, we have $0\geq {\cal L}\tilde G$. Combining the lower bound (\ref{lower bound}) for ${\cal L}\tilde G$ with the preceding inequalities and dropping the second and last terms, which are non-negative, we obtain
\bea
0\geq
{\theta^2\over 10 nP}{\cal F}\lambda_1^2-C_{12}{\cal F}
-{3\over 2}F^{k\bar k}{|\tilde\lambda_{1,\bar k}|^2\over\lambda_1^2}
+
\phi''F^{q\bar q}|\na_{\bar q}|\na u|^2|^2
+\varphi'(\tilde v)(F^{k\bar k}\na_k\na_{\bar k}\tilde v-\p_t\tilde v),
\nonumber\\
\eea
where $C_{12}=C_9+C_{10}+C_{11}$, depending on $\|\chi\|_{C^2},\|\alpha\|_{C^2}, \| \psi\|_{C^2}, \|\p_t u \|_{L^\infty}$  and $n$. Since we are at a critical point of $\tilde G$, we also have $\na \tilde G=0$, 
and hence
\bea
{\tilde \lambda_{1,\bar k}\over \lambda_1}+\phi' \na_{\bar k}|\na u|^2
+
\varphi' \p_{\bar k}\tilde v=0
\eea
which implies
\bea
{3\over 2}F^{k\bar k}|{\tilde \lambda_{1,\bar k}\over \lambda_1}|^2
&=&{3\over 2}F^{k\bar k}|\phi' \na_{\bar k}|\na u|^2
+
\varphi' \p_{\bar k}\tilde v|^2
\leq
2F^{k\bar k}(\phi')^2|\na_{\bar k}|\na u|^2|^2
+
4F^{k\bar k}(\varphi')^2|\na_{\bar k}\tilde v|^2
\nonumber\\
&\leq &
F^{k\bar k}\phi''|\na_{\bar k}|\na u|^2|^2+C_{13}{\cal F}P,
\eea
where $C_{13}$ depending on $\| \tilde{v} \|_{L^\infty}$ and $ \|\na \U \|_{L^\infty}$. 
Since $\varphi'(\tilde v)$ is bounded in terms of $\| \tilde{v} \|_{L^\infty}$ and $ \|\na \U \|_{L^\infty}$, and $|F^{k\bar k}\na_k\na_{\bar k}\tilde v-\p_t\tilde v|\leq C_{14}{\cal F}\lambda_1+C_{13}$, where $ C_{14}$ depending on $ \| \p_t v\|_{L^\infty}$ and $\|\p\bar \p \U\|_{L^\infty}$, we arrive at
\bea
0\geq
{\theta^2\over 10 nP}{\cal F}\lambda_1^2
-C_{15}P{\cal F},
\eea
where $C_{15}$ depends on $ \|\chi\|_{C^2},\|\alpha\|_{C^2}, n,\|\psi\|_{C^2}, \|\p\bar \p \U\|_{L^\infty}, \|\na \U \|_{L^\infty}, \| \tilde{v} \|_{L^\infty}, \| \p_t v\|_{L^\infty}$ and $\|\p_t u\|_{L^\infty}$.
This  implies the desired estimate $\lambda_1\leq \tilde{C} \,P$.

\subsubsection
{The key estimate provided by subsolutions}

In the second case when $\theta\lambda_1>-\lambda_n$, we need to use the following key property of subsolutions.

\begin{lemma}
\label{Festimate}
Let $\U$ be a subsolution of the equation (\ref{parabolic eq}) in the sense of Definition \ref{subsolution} with the pair $(\delta,K)$. Then there exists a constant $C=C(\delta,K)$, so that, if
$|\lambda[u]-\lambda[\U]|>K$ with $K$ in Definition \ref{subsolution}, then either
\bea
\label{ineq_C1}
F^{pq}(A[u])
(A^p{}_{q}[\U]-A^p{}_{q}[u])-(\p_t\U-\p_tu)>C\, {\cal F}
\eea
or we have for any $1\leq i\leq n$,
\bea
F^{ii}(A[u])>C\,{\cal F}.
\eea
\end{lemma}
{\it Proof.} The proof is an adaptation of the one for
the elliptic version \cite[Proposition 6]{Sze}(see also \cite{Gb2} for a similar argument). However, because of the time parameter $t$ which may tend to $\infty$, we need to produce explicit bounds which are independent of $t$.
As in \cite{Sze},  it suffices to prove that
\bea
\sum_{i=1}^{n}f_i(\lambda[u]) (\lambda_i[\underline{u}]-\lambda _i[u]) -(\p_t\U-\p_t u) >C\F.
\eea

For any $(z_0,t_0)\in X\times [0,T'] $, since $\underline{u}$ is a $C$-subsolution as in Definition \ref{subsolution}, the set
$$A_{z_0,t_0} =\{ (w, s)| \, w + \frac{\delta}{2}  I\in \overline  \Gamma_n, s\geq -\delta,\,  \,  f(\lambda[\underline{u} (z_0,t_0)] + w) -\p_t\U (z_0,t_0)+s \leq \psi(z_0) \} $$
is compact, and $A_{z_0,t_0} \subset  B_{n+1}(0,K).$
  For any $(w,s)\in A_{z_0,t_0}$,   then the set
$$
C_{w,s}=\{ v \in \R^{n}| \exists  r >0, \, w+ r v\in -\delta I+\Gamma_n, \, \, \,  f(\lambda[\underline{u}(z_0,t_0)] +w+rv) -\p_t\U(z_0,t_0) +s =\psi(z_0) \}
$$ is  a cone with vertex at the origin. 
 
 \medskip
  We claim that $C_{w,s}$ is stricly larger than $\Gamma_n$. Indeed, for any $v\in \Gamma_n$,  we can choose $r>0$ large enough so that $ |w+rv|>K $, then by the definition of $C$-subsolution,   at $(z_0,t_0)$
  $$f(\lambda[\underline u]+w+rv ) -\p_t\U+s> \psi (z_0).$$ Therefore there exist $r' >0$ such that $f(\lambda[\underline u]) +w+r' v ) -\p_t\U+s=\psi (z_0)$, hence $v\in C_{w,s}$. This implies that $\Gamma_n\subset C_{w,s} $.  Now, for any pair $(i,j)$ with $i\neq j$ and $i,j=1,\ldots,n$, we choose
 $v^{(i,j)}:=(v_1,\ldots,  v_n)$ with $v_i=K+\delta $  and $v_j=-\delta/3$ and $v_k=0$ for $k\neq i,j$,  then we have $w+ v^{(i,j)}\in  -\delta{\bf 1} +\Gamma_n$. By the definition of $ C$-subsolution, we also have, at $(z_0,t_0)$
 $$f(\lambda[\underline u]) +w+ v^{(i,j)} ) -\p_t\U+s> \psi (z_0),$$
hence  $v^{(i,j)}\in C_{w, s}$ for any pair $(i,j)$. 
  
 \medskip
 Denote by $C^*_{w,s}$ the dual cone of $C_{w,s}$,

$$C_{w,s}^*=
\{ x\in \R^n: \langle x,y\rangle >0, \,\forall y\in C_{w,s}  \}.$$
We now prove that  that there is an $\e>0$ such that if $x=(x_1,\ldots,x_n)\in C^*_{w,s}$ is a unit vector, then $x_i>\e,\,\forall i=1,\ldots n$.  First we remark that $x_i>0,\forall i=1,n$ since $\Gamma_n\subset  C_{w,s}$ Suppose that $x_1$ is the smallest element between $x_i$, then  $\langle x,v^{(1,j) }\rangle >0,$ implies that $(K+\delta)x_1\geq { \delta\over 3} x_j $, hence $ (K+\delta)^2x_1^2\geq (\delta^2/9) x^2_j, \forall j=2,\ldots,n $, so $n (K+\delta)^2x_1^2\geq \delta^2/9$. Therefore we can choose $\e={\delta^2\over 9n(K+\delta)^2}$.

\medskip
\noindent
Fix $(z_1,t_1)\in X\times [0,T']$ such that  at this point $|\lambda[u] -\lambda[\underline{u}]|>K$.  Let $\mathcal{T}$ be the tangent plane to $\{(\lambda,\tau) |\, f(\lambda)+\tau=\sigma\} $ at $(\lambda[u(z_1,t_1)],-\p_t u(z_1,t_1))$. There are two cases:

\smallskip
1) There is some point $(w,s)\in A_{z_1,t_1}$ such that  at $(z_1,t_1)$
$$ (\lambda[\underline{u}] +w, -\p_t\U+s) \in \mathcal{T},$$
i.e 
\begin{equation}
\label{tangent}
\nabla f(\lambda[u]) . (\lambda[\underline{u}] +w-\lambda[u]) +(-\p_t\U+s+\p_t u)=0.
\end{equation}
Now for any  $v\in C_{w,s}$, there exist $r >0$ such that $f(\lambda[\underline{u}] +w+r v) -\p_t\U+s=\psi(z)$, this implies that  $$\nabla f(\lambda[u]) . (\lambda[\underline{u}] +w +r  v -\lambda[u]) +(-\p_t\U+s+\p_t u )>0,$$
so combing with (\ref{tangent}) we get
 $$\nabla f(\lambda[u]) . v> 0. $$
It follows that at $(z_1,t_1)$ we have $\nabla f(\lambda[u])(z,t)\in C^*_{w,s}$, so  $ f_i(\lambda[u])\geq  \e \nabla f(\lambda[u]), \forall i=1,\ldots
,n$, hence 
$$f_i(\lambda[u])> \frac{\e}{\sqrt{n}} \sum_p f_p(\lambda[u]),\forall i=1,\ldots,n,$$
where $$\e={\delta^2\over 9n(K+\delta)^2}.$$
\smallskip
2) Otherwise, we observe that if  $ A_{z_1,t_1} \neq \emptyset$, then $(w_0,s_0)=(-\delta/2 ,\ldots,-\delta/2, -\delta )\in A_{z_1,t_1} $ and at $(z_1,t_1)$, $(\lambda[\underline{u}] -w_0,-\underline u_t +s_0 )$ must lie above $\mathcal{T}$ in the sense that  
\begin{equation}
\label{tangent}
 (\nabla f (\lambda[u]),1).(\lambda[\underline u]+w_0-\lambda[u],-\p_t\U+s_0+ \p_t u)>0, \, {\rm at}\, (z_1,t_1).
\end{equation}
 Indeed,  if it is not the case, using the monotonicity of $f$ we can find $v\in \Gamma_n $ such that  $(\lambda[\U]+w_0+v, -\p_t\U+s_0 )\in \mathcal{T}$, so the concavity of $(\lambda,\tau)\mapsto f(\lambda)+\tau$ implies that $(w_0+v,s_0)$ is in $ A_{z_1,t_1}$ and then  satisfies the first case, this gives a contradiction. Now it follows from (\ref{tangent}) that at $(z_1,t_1)$
 \begin{eqnarray*}
  (\nabla f (\lambda[u]),1).(\lambda[\underline u]-\lambda[u],-\p_t\U+  \p_t u)&\geq&-\nabla f (\lambda[u]).w_0-s_0\\
  &=&(\delta/2) \F + \delta\geq (\delta/2)\F,
 \end{eqnarray*}
where $\F=\sum_i f_i(\lambda[u])>0$. This means
\begin{equation}
\label{ineq_C2}
\sum_{i=1}^{n}f_i(\lambda[u]) (\lambda[\underline{u}]-\lambda [u]) -(\p_t\U-u_t) >( \delta/2) \F
\end{equation}
as required. 

\medskip Now if $A_{z_1,t_1}=\emptyset$, then at  $(z_1,t_1)$ $$f(\lambda[\U] +w_0)-\p_t \U +s_0>\psi(z_1),$$ 
hence we also have that $(\lambda[\underline{u}] +w_0,-\p_t \U +s_0 )$  lies  above $\mathcal{T}$ using the concavity of $(\lambda,\tau)\mapsto f(\lambda)+\tau$. By the same argument above, we also obtain the inequality (\ref{ineq_C2}). 

\smallskip

So we get the desired inequalities. Q.E.D.

\subsubsection{Case 2: $\theta\lambda_1>-\lambda_n$}

Set
\bea
I=\{i;\ F^{i\bar i}\geq \theta^{-1}F^{1\bar 1}\}.
\eea
At the maximum point $\p_{\bar k}\tilde G=0$, and we can write
\bea
-\sum_{k\not\in I}F^{k\bar k}{|\tilde \lambda_{1,\bar k}|^2\over \lambda_1^2}
&=&
-\sum_{k\notin I}F^{k\bar k}|\phi'\na_{\bar k}|\na u|^2+\varphi'\p_{\bar k}\tilde v|^2
\nonumber\\
&\geq&
-2(\phi')^2\sum_{k\notin I}F^{k\bar k}|\na_{\bar k}|\na u|^2|^2
-2(\varphi')^2\sum_{k\notin I}F^{k\bar k}|\na_{\bar k}\tilde v|^2
\nonumber\\
&\geq&
-\phi''\sum_{k\notin I}F^{k\bar k}|\na_{\bar k}|\na u|^2|^2
-2(\varphi')^2\theta^{-1}F^{1\bar 1}P-C_{16}{\cal F},
\eea
where $C_{16}$ depends on $\|\na \U \|_{L^\infty} $ and $\| \tilde{v}\|_{L^\infty}$. 
On the other hand,
\bea
-2\theta\sum_{k\in I}F^{k\bar k}{|\tilde \lambda_{1,\bar k}|^2\over \lambda_1^2}
\geq
-2\theta\phi''\sum_{k\in I}F^{k\bar k}|\na_{\bar k}|\na u|^2|^2
-
4\theta(\varphi')^2\sum_{k\in I}F^{k\bar k}|\na_{\bar k}\tilde v|^2.
\eea
Choose $0<\theta<<1$ such that $4\theta(\varphi')^2\leq {1\over 2}\varphi''$. Then (\ref{lower bound}) implies that
\bea
0&\geq&
-{1\over\lambda_1}F^{l\bar k,s\bar r}\na_{\bar 1}g_{\bar kl}\na_1g_{\bar rs}
-
(1-2\theta)\sum_{k\in I}F^{k\bar k}{|\tilde \lambda_{1,\bar k}|^2\over \lambda_1^2}
\nonumber\\
&&
-C{1\over\lambda_1}F^{k\bar k}|\tilde\lambda_{1,\bar k}|
+
{1\over 8P}F^{q\bar q}(|\na_q\na u|^2+|\na_q\bar\na u|^2)
\nonumber\\
&&
+{1\over 2}\varphi'' F^{k\bar k}|\na_{\bar k}\tilde v|^2
+
\varphi'(F^{k\bar k}\na_k\na_{\bar k}\tilde v-\p_t\tilde v)
-2(\varphi')^2\theta^{-1}F^{1\bar 1}P-C_{17}{\cal F},
\eea
where $C_{17}$ depend on $ \|\chi\|_{C^2},\|\alpha\|_{C^2},n,\| \psi \|_{C^2}, \|\p_t u\|_{L^\infty}, \|\tilde{v}\|_{L^\infty}$  and $\|\na \U \|_{L^\infty} $.
The concavity of $F$ implies that
\bea
F^{l\bar k,s\bar r}\na_{\bar 1}g_{\bar kl}\na_1g_{\bar rs}
\leq
\sum_{k\in I}{F^{1\bar 1}-F^{k\bar k}\over\lambda_1-\lambda_k}
|\na_1 g_{\bar 1k}|^2
\eea
since ${F^{1\bar 1}-F^{k\bar k}\over\lambda_1-\lambda_k}\leq 0$. Moreover, for $k\in I$, we have $F^{1\bar 1}\leq \theta F^{k\bar k}$, and the assumption $\theta\lambda_1\geq -\lambda_n$ yields
\bea
{1-\theta\over\lambda_1-\lambda_k}
\geq {1-2\theta\over\lambda_1}.
\eea
It follows that
\bea
\sum_{k\in I}{F^{1\bar 1}-F^{k\bar k}\over\lambda_1-\lambda_k}
|\na_1 g_{\bar 1k}|^2
\leq
-
\sum_{k\in I}{(1-\theta)F^{k\bar k}\over\lambda_1-\lambda_k}
|\na_1 g_{\bar 1k}|^2
\leq
-{1-2\theta\over\lambda_1}\sum_{k\in I}F^{k\bar k}|\na_1 g_{\bar 1k}|^2.
\eea
Combining with the previous inequalities, we obtain
\bea
0&\geq&
-(1-2\theta)
\sum_{k\in I}F^{k\bar k}{|\tilde\lambda_{1,\bar k}|^2-|\na_1g_{\bar 1k}|^2\over\lambda_1^2}-C_{17}{\cal F}
\nonumber\\
&&
-{C_4\over\lambda_1}F^{k\bar k}|\tilde\lambda_{1,\bar k}|
+
{1\over 8P}F^{q\bar q}(|\na_q\na u|^2+|\na_q\bar\na u|^2)
\nonumber\\
&&
+{1\over 2}\varphi''F^{k\bar k}|\na_{\bar k}\tilde v|^2
+
\varphi'(F^{k\bar k}\na_k\na_{\bar k}\tilde v-\p_t\tilde v)
-2(\varphi')^2\theta^{-1}F^{1\bar 1}P.
\eea
Since $ \na_1 g_{\bar 1  k}=\tilde{\lambda}_{1,\lambda}+ O(\lambda_1)$, we have
\begin{equation}
-(1-2\theta)
\sum_{k\in I}F^{k\bar k}{|\tilde\lambda_{1,\bar k}|^2-|\na_1g_{\bar 1k}|^2\over\lambda_1^2}\geq - C_{18}{\cal F}
\end{equation}
where $C_{18}$ depends on $\|\chi\|_{C^2}$ and $\|\alpha\|_{C^2}$.
Next, using again the equations for critical points, we can write
\bea
{C_4\over\lambda_1}F^{k\bar k}|\tilde\lambda_{1,\bar k}|
&=&
{C_4\over\lambda_1}F^{k\bar k}| \phi' \na_{\bar k}| \na u|^2 +\varphi'\na_{\bar k}\tilde v|
\\
&\leq &
{1\over 2K^{1\over 2}}\sum F^{k\bar k}(|\na_{\bar k}\na_p u |+|\na_{\bar k}\na_{\bar p} u|)
+
C_\varepsilon |\varphi'|F^{k\bar k}|\na_{\bar k}\tilde v|^2
+
\varepsilon C_{19}|\varphi'|{\cal F}
+
C_{20}{\cal F},\nonumber
\eea
where $C_{19}$ and $C_{20}$ depend on $C_4$.
Accordingly, the previous inequality implies
\bea
0&\geq&
 {1\over 10K}F^{q\bar q}(|\na_q\na u|^2+|\na_q\bar\na u|^2)
+
{1\over 2}\varphi'' F^{k\bar k}|\na_{\bar k}\tilde v|^2
+
\varphi'(F^{k\bar k}\na_k\na_{\bar k}\tilde v-\p_t\tilde v)
\nonumber\\
&&
-2(\varphi')^2\theta^{-1}F^{1\bar 1}P
-
C_\varepsilon |\varphi'|F^{k\bar k}|\na_{\bar k}\tilde v|^2- \varepsilon C_{19}|\varphi'|{\cal F}-C_{21}{\cal F},
\eea
where $C_{21}$ depending only on $\|\chi\|_{C^2},\|\alpha\|_{C^2},n, \|\psi\|_{C^2},
\|\p_tv\|_{C^0},\|\tilde{v}\|_{L^\infty}$, $\| \p_t u\|_{L^\infty}$ and $\|\na \U\|_{L^\infty}$.
Finally we get
\bea
\label{below2}
0&\geq &
F^{1\bar 1}({\lambda_1^2\over 20 P}-2(\varphi')^2\theta^{-1}P)+
({1\over 2}\varphi'' -C_\varepsilon|\varphi'|)F^{k\bar k}|\na_{\bar k}\tilde v|^2
\nonumber\\
&&
-\varepsilon C_{19}|\varphi'|{\cal F}+\varphi'
(F^{k\bar k}\na_k\na_{\bar k}\tilde v-\p_t\tilde v)-C_{21}{\cal F}.
\eea
We now apply Lemma \ref{Festimate}. Fix $\delta$ and $K$ as in Definition \ref{subsolution}, if $\lambda_1>K$, then there are two possibilities:

\smallskip
$\bullet$  Either $F^{k\bar k}(\U_{\bar kk}-u_{\bar kk})+(\p_tu-\p_t\U)\geq \kappa {\cal F}$, for some $\kappa$ depending only on $\delta$ and $ K$, equivalently,
\bea
F^{k\bar k}\na_k\na_{\bar k}\tilde v-\p_t\tilde v-\int_X \p_tv\alpha^n\leq -\kappa {\cal F}+C_{22}{\cal F},
\eea
where $C_{22}$ depends on $ \| \p_t v \|_{L^\infty}$.
Since $\varphi'<0$, we find
\bea
0&\geq&
F^{1\bar 1}({\lambda_1^2\over 20 P}-2(\varphi')^2\theta^{-1}P)+
({1\over 2}\varphi''-C_\varepsilon|\varphi'|)F^{k\bar k}|\na_{\bar k}\tilde v|^2
\nonumber\\
&&
-C_{23}{\cal F}-\varepsilon C_{19}|\varphi'|{\cal F}
-
\varphi'\kappa{\cal F}
\eea
with $C_{23}$ depending only on $n,\|\chi\|_{C^2},\|\alpha\|_{C^2}, \|\psi\|_{C^2},
\|\p_tv\|_{L^\infty}$, $\|\tilde{v}\|_{L^\infty},\|\p_t u\|_{L^\infty}$ and $\|\na \U\|_{L^\infty}$. We first choose $\varepsilon$ small enough so that $\varepsilon C_{19}<\kappa/2$, then $D_2$ large enough so that $\varphi''>2C_\varepsilon |\varphi'|$. We obtain
\bea
0\geq 
F^{1\bar 1}({\lambda_1^2\over 20 P}-2(\varphi')^2\theta^{-1}P)
-C_{23}{\cal F}-{1\over 2}\varphi'\kappa {\cal F}.
\eea
We now choose $D_1$ large enough (depending on $\|\tilde v\|_{L^\infty}$) so that 
$-C_{23}-{1\over 2}\varphi'\kappa>0$. Then
\bea
{\lambda_1^2\over 20 P}\leq 2(\varphi')^2\theta^{-1}P
\eea
and the desired upper bound for $\lambda_1/P$ follows.

\medskip
$\bullet$ Or $F^{1\bar 1}\geq \kappa {\cal F}$. With $D_1$, $D_2$, and $\theta$ as above, the inequality (\ref{below2}) implies
\bea
0\geq
\kappa{\cal F}({\lambda_1^2\over 20 P}-2(\varphi')^2\theta^{-1}P)
-C_{24}{\cal F}-\varphi'F^{k\bar k}g_{\bar kk},
\eea
with $C_{24}$ depending only on $\|\chi\|_{C^2},\|\alpha\|_{C^2},n, \|\psi\|_{C^2},
\|\p_tv\|_{L^\infty},\|\tilde{v}\|_{L^\infty}, \|\p_t u\|_{L^\infty}$, $\|\na \U\|_{L^\infty}$, and $\|i\p\bar\p\U\|_{L^\infty}$. Since $F^{k\bar k}g_{\bar kk}\leq {\cal F}\lambda_1$, we can divide by ${\cal F}P$ to get
\bea
0\geq
\kappa{\lambda_1^2\over 20 P^2}-C_{25}(1+{1\over P}+{\lambda_1\over P})
\eea
with a constant $C_{25}$ depending only on $\|\chi\|_{C^2},\|\alpha\|_{C^2}, n,\|\psi\|_{C^2},\|\tilde{v}\|_{L^\infty},\|\p_tv\|_{L^\infty},\|\p_t u\|_{L^\infty}$, $\|\na \U\|_{L^\infty}$, and $\|i\p\bar\p\U\|_{L^\infty}$. Thus we obtain the desired bound for $\lambda_1/P$.

\medskip
It was pointed out in \cite{Sze} that, under an extra concavity condition on $f$, $C^2$ estimates can be derived directly from $C^0$ estimates in the elliptic case, using a test function introduced in \cite{PS09}. The same holds in the parabolic case, but we omit a fuller discussion.

\subsection{$C^1$ Estimates}

The $C^1$ estimates are also adapted from \cite{Sze}, which reduce the estimates by a blow-up argument to a key Liouville theorem for Hessian equations due to Sz\'ekelyhidi \cite{Sze} and Dinew and Kolodziej \cite{DK}.

\begin{lemma}\label{C1}
There exist a constant $C>0$, depending on  $\underline{u}$, $\|\p_t u \|_{L^\infty(X\times[0,T))}$, $\| \tilde{u} \|_{L^\infty(X\times[0
,T))} $  $\|\alpha\|_{C^2}, \chi,\psi $  and the constant $\tilde{C}$ in Lemma \ref{C2} such that 
\begin{equation}\label{bound gradient}
\sup_{X\times [0,T)}|\na u|^2_\alpha\leq C.
\end{equation}
\end{lemma}

\noindent
{\it Proof.}
Assume by contradiction that (\ref{bound gradient}) does not hold. Then there exists a sequence $(x_k,t_k)\in X\times [0,T)$ with $t_k\rightarrow T$ such that
$$
\lim_{k\rightarrow\infty} |\nabla u(t_k,x_k)|_\alpha =+\infty.
$$
We can assume further that 
$$R_k=|\nabla  u(x_k,t_k)|_\alpha =\sup_{X\times [0,t_k]}|\nabla  u(x,t)|_\alpha,\quad as \quad k\rightarrow +\infty,$$
and $\lim_{k\rightarrow \infty} x_k=x$. 

\medskip
\noindent
Using localization, we choose a coordinate chart $\{U, (z_1,\ldots,z_n)\}$ centered at $x$, identifying with the ball $B_2(0)\subset \C^n$ of radius 2 centered at the origin such that $\alpha (0)=\beta$, where $\beta=\sum_j idz^j\wedge d\bar z^j$.   We also assume that $k$ is sufficiently large so that $z_k:=z(x_k)\in B_1(0)$.

\medskip 
\noindent
Define the following maps
\begin{eqnarray*}
&&\Phi_k:\C^n\rightarrow \C^n,\quad \Phi_k(z):=R_k^{-1}z+z_k,\\
&& \tilde u_k: B_{R_k}(0)\rightarrow \R, \quad \tilde u_k(z):=\tilde u( \Phi_k(z),t_k)=\tilde u(R^{-1}_kz+z_k,t_k),
\end{eqnarray*}
where $\tilde{u}=u-\int_X u\,\alpha^n$.
Then the equation $$u_t=F(A)-\psi(z),$$ implies that
\begin{equation}\label{new equation}
 f\left(R_k^2 \lambda [\beta^{i\bar p}_k (\chi_{k,\bar p j} +\tilde{u}_{k,\bar p j})]\right) = \psi (R_k^{-1}z+z_k)+u_t(\Phi_k(z),t_k) ,
\end{equation}
where $\beta_k:= R^2_k\Phi^*_k\alpha, \chi_k:=\Phi_k^* \chi$. 
Since $ \beta_k  \rightarrow \beta$, and $\chi_k(z,t)\rightarrow 0$,  in $C^{\infty}_{loc}$  as $k\rightarrow \infty$, we get
\begin{equation}\label{asymptotic 1}
 \lambda [\beta^{i\bar p}_k (\chi_{k,\bar p j} +\tilde{u}_{k,\bar p j})]=\lambda (\tilde u_{k,\bar j i})+ O\left(\frac{|z|}{R^2_k}\right).
\end{equation}
By the construction,  we  have
\begin{equation}
\label{bound_01}
\sup_{B_{R_k}(0)} \tilde{u}_k \leq C, \quad \sup_{B_{R_k}(0)} |\nabla  \tilde{u}_k | \leq C
\end{equation}
where $C$ depending on $\|\tilde u\|_{L^\infty}$,
and $$|\nabla  \tilde{u}_k | (0)=R_k^{-1} | \nabla  u_k|_\alpha (x_k)=1. $$
Thanks to Lemma \ref{C2}, we also have that
\begin{equation}
\label{bound_2}
\sup_{B_{R_k}(0)}| \partial\bar{\partial}\tilde{u}_k |_\beta \leq CR_k^{-2}\sup_{X}|\partial \bar\partial u(.,t_k) |_\alpha\leq C' .
\end{equation}
As the argument in \cite{Sze,TW17}, it follows from (\ref{bound_01}), (\ref{bound_2}), the elliptic estimates for $\Delta$ and the Sobolev embedding that for each given $K\subset\C^n$ compact, $0<\gamma< 1$ and $p>1$, there is a constant $C$  such that 
$$\|\tilde{u}_k\|_{C^{1,\gamma}(K)} +\|\tilde{u}_k \| _{W^{2,p} (K)}\leq C.$$
Therefore there is a subsequence of $ \tilde{u}_k$ converges strongly in $C^{1,\gamma}_{loc}(\C^n)$, and weakly in $W^{2,p}_{loc}(\C^n)$ to a function $v$ with $\sup_{\C^n}(|v|+|\nabla  v| )\leq C$ and $\nabla  v(0)\neq 0$, in particular $v$ is not constant. 

The proof can now be completed exactly as in \cite{Sze}. The function $v$ is shown to be a $\Gamma$-solution in the sense of Sz\'ekelyhidi \cite[Definition 15]{Sze}, and the fact that $v$ is not constant contradicts Szekelyhidi's Liouville theorem for $\Gamma$-solutions \cite[Theorem 20]{Sze}, which is itself based on the Liouville theorem of Dinew and Kolodziej \cite{DK}. Q.E.D.

\subsection{Higher Order Estimates}

Under the conditions on $f(\lambda)$, the uniform parabolicity of the equation (\ref{parabolic eq}) will follow once we have established an a priori estimate on $\|i\p\bar\p u\|_{L^\infty}$ and hence an upper bound for the eigenvalues $\lambda[u]$. However, we shall often not have uniform control of $\|u(\cdot,t)\|_{L^\infty}$. Thus we shall require the following version of the Evans-Krylov theorem for uniformly parabolic and concave equations, with the precise dependence of constants spelled out, and which can be proved using the arguments of
Trudinger \cite{Tr}, and more particularly Tosatti-Weinkove \cite{TW10b} and Gill \cite{G}.

\begin{lemma}
\label{C2alpha}
Assume that u is a solution of the equation
(\ref{parabolic eq}) on $X\times [0, T)$ and that there exists a constant
$C_0$ with $\|i\p\bar\p u\|_{L^\infty}\leq C_0$.
Then there exist positive constants $C$ and $\gamma\in (0, 1)$ depending 
only on $\alpha$, $\chi$, $C_0$ and $\|\psi\|_{C^2}$ 
such that
\bea
\|i\p\bar\p u\|_{C^\gamma(X\times[0,T))}\leq C.
\eea
\end{lemma}

Once the $C^\gamma$ estimate for $i\p\bar\p u$ has been established, it is well known that a priori estimates of arbitrary order follow by bootstrap, as shown in detail for the Monge-Amp\`ere equation in Yau \cite{Y}. We omit reproducing the proofs.

\section{Proof of Theorems 1 and 2}
\setcounter{equation}{0}

We begin with the following simple lemma, which follows immediately by differentiating the equation (\ref{parabolic eq}) with respect to $t$, and applying the maximum principle, which shows that the solution of a linear heat equation at any time can be controlled by its initial value:

\begin{lemma}
\label{pt}
Let $u(z,t)$ be a smooth solution of the flow (\ref{parabolic eq}) on any time interval $[0,T)$. Then $\p_tu$ satisfies the following linear heat equation
\bea
\p_t(\p_t u)=F^j{}_k\alpha^{k\bar m}\p_j\p_{\bar m}(\p_tu)
\eea
and we have the following estimate for any $t\in [0,T)$,
\bea
{\rm min}_X (F(A[u_0])-\psi)
\leq \p_tu(t,\cdot)
\leq
{\rm max}_XF(A[u_0]-\psi)
\eea
\end{lemma}

\medskip
\noindent
We can now prove a lemma which provides general sufficient conditions for the convergence of the flow:

\begin{lemma}
\label{convergence lemma}
Consider the flow (\ref{parabolic eq}). Assume that the equation admits a parabolic $C$-subsolution $\underline{u}\in C^{2,1}(X\times [0,\infty))$, and that there exists a constant $C$ independent of time so that
\bea
\label{osc}
{\rm osc}_X u(t,\cdot)\leq C.
\eea
Then a smooth solution $u(z,t)$ exists for all time, and its normalization $\tilde u$ converges in $C^\infty$ to a solution $u_\infty$ of the equation (\ref{equation})
for some constant $c$.

\medskip
\noindent In particular, if we assume further that $\|u\|_{L^\infty( X\times [0,\infty))}\leq C$ and for each $t> 0$, there exists $y=y(t)\in X$  such that $\partial_t u(y,t)=0$, then $u$ converges in $C^\infty$ to a solution $u_\infty$ of the equation (\ref{equation}) 
for the constant $c=0$. 
\end{lemma}

\noindent
{\it Proof of Lemma \ref{convergence lemma}}.
We begin by establishing the existence of the solution for all time. For any fixed $T>0$, Lemma \ref{pt} shows that $|\p_tu|$ is uniformly bounded by a constant $C$. Integrating between $0$ and $T$, we deduce that $|u|$ is uniformly bounded by $C\,T$. We can now apply Lemma \ref{C1}, \ref{C2}, \ref{C2alpha}, to conclude that the function $u$ is uniformly bounded in $C^k$ norm (by constants depending on $k$ and $T$) for arbitrary $k$. This implies that the solution can be extended beyond $T$, and since $T$ is arbitrary, that it exists for all time.

\smallskip

Next, we establish the convergence. For this, we adapt the arguments of Cao \cite{Ca} and especially Gill \cite{G} based on the Harnack inequality.

Since ${\rm osc}_X u(t,\cdot)$ is uniformly bounded by assumption, and since
$\p_tu$ is uniformly bounded in view of Lemma \ref{pt}, we can apply
Lemma \ref{C2} and deduce that the eigenvalues of the matrix $[\chi+i\p\bar\p u]$ are uniformly bounded over the time interval $[0,\infty)$. The uniform ellipticity of the equation
(\ref{heat eq}) follows in turn from the properties (1) and (2) of the function $f(\lambda)$. Next set
\bea
v=\p_t u +A
\eea
for some large constant $A$ so that $v>0$. The function $v$ satisfies the same heat equation
\bea
\label{heat eq}
\p_tv=F^{i\bar j}\p_i\p_{\bar j}v.
\eea
Since the equation (\ref{heat eq}) is uniformly elliptic, by the differential Harnack inequality proved originally in the Riemannian case by Li and Yau in \cite{LY}, and extended to the Hermitian case by Gill \cite{G}, section 6,
it follows that there exist positive constants $C_1, C_2,C_3$, depending only on ellipticity bounds, so that for all $0<t_1<t_2$, we have
\bea
{\rm sup}_X v(\cdot,t_1)
\leq
{\rm inf}_X v(\cdot,t_2)
\left({t_2\over t_1}\right)^{C_2}{\rm exp}
\left({C_3\over t_2-t_1}+C_1(t_2-t_1)\right).
\eea
The same argument as in Cao \cite{Ca}, section 2, and Gill \cite{G}, section 7, shows that this estimate implies the existence of constants $C_4$ and $\eta>0$ so that
\bea \label{osc_v}
{\rm osc}_Xv(\cdot,t)
\leq
C_4 e^{-\eta t}
\eea
If we set
\bea
\tilde v(z,t)=v(z,t)-{1\over V}\int_X v\,\alpha^n=\p_t u(z,t)-{1\over V}\int_X \p_tu\,\alpha^n
=\p_t\tilde u,
\eea
it follows that
\bea 
|\tilde v(z,t)|\leq C_4 e^{-\eta t}
\eea
for all $z\in X$. In particular, 
\bea
\p_t(\tilde u+{C_4\over\eta}e^{-\eta t})
=
\tilde v-C_4 e^{-\eta t}\leq 0,
\eea
and the function $\tilde u(z,t)+{C_4\over\eta}e^{-\eta t}$ is decreasing in $t$. By the assumption (\ref{osc}), this function is uniformly bounded. Thus it converges to a function $u_\infty(z)$. By the higher order estimates in section \S 2, the derivatives to any order of $\tilde u$ are uniformly bounded, so the convergence of $\tilde u+{C_4\over\eta}e^{-\eta t}$ is actually in $C^\infty$ The function $\tilde u(z,t)$ will also converge in $C^\infty$, to the same limit $u_\infty(z)$. Now the function $\tilde u(z,t)$ satisfies the following flow,
\bea
\p_t\tilde u=F(A[\tilde u])-\psi(z)-{1\over V}\int_X \p_tu\,\alpha^n.
\eea
Taking limits, we obtain
\bea
0=F(A[\tilde u_\infty])-\psi(z)-{\rm lim}_{t\to\infty}\int_X \p_t u\,\alpha^n
\eea
where the existence of the limit of the integral on the right hand side follows from the equation. Define the constant $c$ as the value of this limit. This implies the first statement in  Lemma \ref{convergence lemma}.

\medskip Now  we assume that $\|u\|_{L^\infty( X\times [0,\infty))}\leq C$ and for each $t\geq 0$, there exists $y=y(t)\in X$  such that $\partial_t u(y,t)=0$.  By the same argument above,  we have
\bea
{\rm osc}_X \partial_t u (\cdot,t)
\leq
C_4 e^{-\eta t},
\eea
for some $C_4, \eta>0$.  Since for each $t\geq 0$, there exists $y=y(t)\in X$  such that $\partial_t u(y,t)=0$, we imply that for any $z\in X$, 
\bea
|\p_t u(z,t)|=| \p_t u(z,t)-\p_t u(y,t )| \leq {\rm osc}_X \partial_t u (\cdot,t)
\leq
C_4 e^{-\eta t}.
\eea
Therefore by the same argument above, the function $u(z,t)+ \frac{C_4}{\eta} e^{-\eta t} $ converges in $C^\infty$ and $\p_t u$ converges to $0$ as $t\rightarrow +\infty$. 
We thus infer that $u$ converges in $C^\infty $, to $u_\infty$ satisfying the equation 
\bea
F(A[\tilde u_\infty])=\psi(z).
\eea
 Lemma \ref{convergence lemma} is proved.

\bigskip
\noindent
{\it Proof of Theorem \ref{theorem B}}.  Since $f$ is unbounded, the function $\underline{u}=u_0$ is a $C$-subsolution of the flow.  In view of Lemma \ref{convergence lemma}, it suffices to establish a uniform bound for ${\rm osc}_X u(t,\cdot)$. But the flow can be re-expressed as the elliptic equation
\bea \label{elliptic_eq_unbounded}
F(A)=\psi+\p_tu
\eea
where the right hand side $\psi+\p_tu$ is bounded uniformly in $t$, since we have seen that $\p_tu$ is uniformly bounded in $t$. Furthermore, because $f$ is unbounded, the function $\underline{u}=u_0$ is a $C$-subsolution of (\ref{elliptic_eq_unbounded}). By the $C^0$ estimate of \cite{Sze}, the oscillation ${\rm osc}_X u(t,\cdot)$ can be bounded for each $t$ by the $C^0$ norm of the right hand side, and is hence uniformly bounded. Q.E.D.

\bigskip
\noindent
{\it Proof of Theorem \ref{theorem C}}. Again, it suffices to establish a uniform bound in $t$ for ${\rm osc}_X u(t,\cdot)$.

Consider first the case (a). In view of Lemma \ref{pt} and the hypothesis, we have
\bea
\label{u-underline u}
\p_t\underline{u}\geq \p_t u
\eea
on all of $X\times [0,\infty)$. But if we rewrite the flow (\ref{parabolic eq}) as
\bea
\label{elliptic version}
F(A)=\psi+\p_t u
\eea
we see that the condition that $\underline{u}$ be a parabolic $C$-subsolution for the equation (\ref{parabolic eq}) together with (\ref{u-underline u}) implies that $\underline{u}$ is a $C$-subsolution for the equation (\ref{elliptic version}) in the elliptic sense.. We can then apply Sz\'ekelyhidi's $C^0$ estimate for the elliptic equation to obtain a uniform bound for ${\rm osc}_X u(t,\cdot)$.

Next, we consider the case (b). In this case, the existence of a function $h(t)$ with the indicated properties allows us to apply Lemma \ref{C0}, and obtain immediately a lower bound,
\bea
u-\underline{u}-h(t)\geq -C
\eea
for some constant $C$ independent of time. The inequality (\ref{harnack_ineq}) implies than a uniform bound for ${\rm osc}_X u$.

\section{Applications to Geometric Flows}
\setcounter{equation}{0}

Theorems \ref{theorem B} and \ref{theorem C} can be applied to many geometric flows. We should stress that they don't provide a completely independent approach, as they themselves are built on many techniques that had been developed to study these flows. Nevertheless, they may provide an attractive uniform approach.

\subsection{A criterion for subsolutions}

In practice, it is easier to verify that a given function $\U$ on $X\times [0,\infty)$ is a $C$-subsolution of the equation (\ref{parabolic eq}) using the following lemma rather than the original Definition \ref{subsolution}:

\begin{lemma}
\label{criterion}
Let $\U$ be a $C^{2,1}$ admissible function on $X\times [0,\infty)$, with $\|\U\|_{C^{2,1}(X\times [0,\infty))}<\infty$. Then $\U$ is a parabolic $C$-subsolution in the sense of Definition \ref{subsolution} if and only if there exists a constant $\tilde\delta>0$
independent from $(z,t)$ so that
\bea
\label{def_2}
{\rm lim}_{\mu\to +\infty}
f(\lambda[\U(z,t)]+\mu e_i)-\p_t\U(z,t)>\tilde\delta+\psi(z)
\eea
for each $1\leq i\leq n$. 
In particular, if $\underline{u}$ is independent of $t$, then $\U$ is a parabolic $C$-subsolution if and only if
\bea
\label{def_3}
{\rm lim}_{\mu\to +\infty}
f(\lambda[\U(z,t)]+\mu e_i)>\psi(z).
\eea
\end{lemma}

Note that there is a similar lemma in the case of subsolutions for elliptic equations (see \cite{Sze}, Remark 8). Here the argument has to be more careful, not just because of the additional time parameter $t$, but also because the time interval $[0,\infty)$ is not bounded, invalidating certain compactness arguments.

\medskip
\noindent
{\it Proof of Lemma \ref{criterion}.} 
We show first that the condition (\ref{def_2}) implies that $\U$ is a $C$-subsolution. 

\smallskip

We begin by showing that the condition (\ref{def_2}) implies that there exists $\epsilon_0>0$ and $M>0$, so that for all $\epsilon\leq\epsilon_0$, all $\nu>M$, all $(z,t)$, and all $1\leq i\leq n$, we have
\bea
\label{sub1}
f(\lambda[\U(z,t)]-\epsilon I+\nu e_i)
-\p_t\U(z,t)>{\tilde\delta\over 4}+\psi(z).
\eea
This is because the condition (\ref{def_2}) is equivalent to
\bea
f_\infty(\lambda'[\U(z,t)])-\p_t\U(z,t)>\tilde\delta+\psi(z).
\eea
Now the concavity of $f(\lambda)$ implies the concavity of its limit $f_\infty(\lambda')$ and hence the continuity of $f_\infty(\lambda')$. Furthermore,
the set 
\bea
\Lambda=\overline{\{\lambda [\underline u(z,t)], \forall (z,t)\in X\times [0,\infty)\}},
\eea 
as well as any of its translates by $-\epsilon I$ for a fixed $\epsilon$ small enough, is compact in $\Gamma$. So are their projections on ${\bf R}^{n-1}$. By the uniform continuity of continuous functions on compact sets, it follows that there exists $\epsilon_0>0$ so that
\bea
\label{sub2}
f_\infty(\lambda'[\U(z,t)]-\epsilon I)-\p_t\U(z,t)>{\tilde\delta\over 2}+\psi(z)
\eea
for all $(z,t)$ and all $\epsilon\leq \epsilon_0$. But $f_\infty$ is the continuous limit of a sequence of monotone increasing continuous functions
\bea
f_\infty(\lambda'-\epsilon I)={\rm lim}_{\nu\to\infty}
f(\lambda-\epsilon I+\nu e_i).
\eea
By Dini's theorem, the convergence is uniform over any compact subset. Thus there exists $M>0$ large enough so that $\nu>M$ implies that
\bea
\label{sub3}
f(\lambda[\U(z,t)]-\epsilon I+\nu e_i)>f_\infty(\lambda'[\U(z,t)]-\epsilon I)
-{\tilde\delta\over 4}
\eea
for all $(z,t)$ and all $\epsilon\leq\epsilon_0$. The desired inequality (\ref{sub1}) follows from (\ref{sub2}) and (\ref{sub3}).

\smallskip
Assume now that $\U$ is not a $C$-subsolution. Then there exists $\epsilon_m$, $\nu_m$, $\tau_m$, with $\epsilon_m\to 0$, $\nu_m\in -\epsilon_m I+\Gamma_n$, $\tau_m>-\epsilon_m$, and $|\tau_m|+|\nu_m|\to\infty$, so that
\bea
\label{sub4}
f(\lambda[\U(z_m,t_m)]
+\nu_m)-\p_t\U(z_m,t_m)+\tau_m
=\psi(z_m,t_m).
\eea
Set $\nu_m=-\epsilon_m +\mu_m$, with $\mu_m\in \Gamma_n$. Then we can write
\bea
\tau_m
&=&
-f(\lambda[\U(z_m,t_m)]-\epsilon_mI
+\mu_m)
+\p_t\U(z_m,t_m)+\psi(z_m,t_m)
\nonumber\\
&\leq&
-f(\lambda[\U(z_m,t_m)]-\epsilon_mI)
+\p_t\U(z_m,t_m)+\psi(z_m,t_m)
\eea
which is bounded by a constant. Thus we must have $|\nu_m|$ tending to $+\infty$, or equivalently, $|\mu_m|$ tending to $+\infty$.

\smallskip
By going to a subsequence, we may assume that there is an index $i$ for which the $i$-th components $\mu_m^i$ of the vector $\mu_m$ tend to $\infty$ as $m\to\infty$. By the monotonicity of $f$ in each component, we have
\bea
f(\lambda[\U(z_m,t_m)]-\epsilon_m I+\mu_m^ie_i)
-\p_t\U(z_m,t_m)
&\leq&
f(\lambda[\U(z_m,t_m)]-\epsilon_m I+\mu_m)
-\p_t\U(z_m,t_m)
\nonumber\\
&=&
f(\lambda[\U(z_m,t_m)]+\nu_m)-\p_t\U(z_m,t_m).
\nonumber
\eea
In view of (\ref{sub1}), the left hand side is $\geq {\tilde\delta\over 4}+\psi(z_m,t_m)$ for $\mu_m^i$ large and $\epsilon_m$ small enough. On the other hand, the equation (\ref{sub4}) implies that the right hand side is equal to $\psi(z_m,t_m)-\tau_m$. Thus we obtain
\bea
{\tilde\delta\over 4}+\psi(z_m,t_m)
\leq
\psi(z_m,t_m)-\tau_m
\leq
\psi(z_m,t_m)+\epsilon_m.
\eea
Hence ${\tilde\delta\over 4}\leq\epsilon_m$, which is a contradiction, since $\epsilon_m\to 0$.

\smallskip
Finally, we show that if $\U$ is a subsolution, it must satisfy the condition (\ref{def_2}). Assume otherwise. Then there exists an index $i$ and a sequence $\delta_m\to 0$ and points $(z_m,t_m)$ so that
\bea
{\rm lim}_{\nu\to\infty}
f(\lambda[\U(z_m,t_m)]+\nu e_i)
-\p_t\U(z_m,t_m)\leq
\delta_m+\psi(z_m).
\eea
Since $f$ is increasing in $\nu$, this implies that for any $\nu\in{\bf R}_+$, we have
\bea
f(\lambda[\U(z_m,t_m)]+\nu e_i)
-\p_t\U(z_m,t_m)\leq
\delta_m+\psi(z_m).
\eea
For each $\nu\in {\bf R}_+$, define $\tau_m$ by the equation
\bea
f(\lambda[\U(z_m,t_m)]+\nu e_i)
-\p_t\U(z_m,t_m)+\tau_m=\psi_m.
\eea
The previous inequality means that $\tau_m\geq -\delta_m$, and thus the pair
$(\tau_m, \mu=\nu e_i)$ satisfy the equation (\ref{sub0}). Since we can take $\nu\to +\infty$, this contradicts the defining property of $C$-subsolutions. The proof of Lemma \ref{criterion} is complete.

\subsection{Sz\'ekelyhidi's theorem}

Theorem \ref{theorem C} can be applied to provide a proof by parabolic methods of the following theorem originally proved by Sz\'ekelyhidi \cite{Sze}:

\begin{corollary}
\label{Cor1}
Let $(X,\alpha)$ be a compact Hermitian manifold, and $f(\lambda)$ be a function satisfying the conditions (1-3) spelled out in \S 1 and in the bounded case. Let $\psi$ be a smooth function on $X$. If there exists an admissible function $u_0$ with $F(A[u_0])\leq \psi$, and if the equation $F(A[u])=\psi$ admits a $C$-subsolution in the sense of \cite{Sze}, then the equation $F(A[u])=\psi+c$ admits a smooth solution for some constant $c$.
\end{corollary}

\smallskip
\noindent
{\it Proof of Corollary \ref{Cor1}}. It follows from Lemma \ref{criterion} that a $C$-subsolution in the sense of \cite{Sze} of the elliptic equation $F(A[u])=\psi$ can be viewed as a time-independent parabolic $C$-subsolution $\underline{u}$ of the equation (\ref{parabolic eq}). Consider this flow with initial value $u_0$. Then
\bea
\p_t\underline{u}=0
\geq F(A[u_0])-\psi.
\eea
Thus condition (a) of Theorem \ref{theorem C} is satisfied, and the corollary follows.

\subsection{The K\"ahler-Ricci flow and the Chern-Ricci flow}

On K\"ahler manifolds $(X,\alpha)$ with $c_1(X)=0$,
the K\"ahler-Ricci flow is the flow $\dot g_{\bar kj}=-R_{\bar kj}$. For initial data in the K\"ahler class $[\alpha]$, the evolving metric can be expressed as
$g_{\bar kj}=\alpha_{\bar kj}+\p_j\p_{\bar k}\varphi$, and the flow is equivalent to the following Monge-Amp\`ere flow,
\bea
\label{MA}
\p_t\varphi=\log{(\alpha+i\p\bar\p\varphi)^n\over\alpha^n}-\psi(z)
\eea
for a suitable function $\psi(z)$ satisfying the compatibility condition $\int_X e^\psi\alpha^n=\int_X\alpha^n$. The convergence of this flow was proved by Cao\cite{Ca}, thus giving a parabolic proof of Yau's solution of the Calabi conjecture \cite{Y}. We can readily derive Cao's result from Theorem \ref{theorem B}:

\begin{corollary}
\label{Cor2}
For any initial data, the normalization $\tilde{\varphi}$ of the flow (\ref{MA}) converges in $C^\infty$ to a solution of the equation $(\alpha+i\p\bar\p\varphi)^n=e^\psi\alpha^n$.
\end{corollary}

\noindent
{\it Proof of Corollary \ref{Cor2}}. The Monge-Amp\`ere flow (\ref{MA}) corresponds to the equation (\ref{parabolic eq}) with $\chi=\alpha$, $f(\lambda)=\log \prod_{j=1}^n\lambda_j$, and $\Gamma$ being the full octant $\Gamma_n$. It is straightforward that $f$ satisfies the condition (1-3) in \S 1. In particular $f$ is in the unbounded case, and Theorem \ref{theorem B} applies, giving the convergence of the normalizations $\tilde u(\cdot,t)$ to a smooth solution of the equation 
$(\alpha+i\p\bar\p\varphi)^n=e^{\psi+c}\alpha^n$ for some constant $c$. Integrating both sides of this equation and using the compatibility condition on $\psi$, we find that $c=0$. The corollary is proved.

\medskip
The generalization of the flow (\ref{MA}) to the more general set-up of a compact
Hermitian manifold $(X,\alpha)$ was introduced by Gill \cite{G}. It is known as the Chern-Ricci flow, with the Chern-Ricci tensor $Ric^C(\o)=-i\p\bar\p\log \o^n$ playing the role of the Ricci tensor in the K\"ahler-Ricci flow (we refer to \cite{TW13,TW15,TW15b,To2} and references therein). Gill proved the convergence of this flow, thus providing an alternative proof of the generalization of Yau's theorem proved earlier by Tosatti and Weinkove \cite{TW1}. Generalizations of Yau's theorem had attracted a lot of attention, and many partial results had been obtained before, including those of Cherrier \cite{Ch}, Guan-Li \cite{GL}, and others. Theorem \ref{theorem B} gives immediately another proof of Gill's theorem:

\begin{corollary}
\label{Cor3}
For any initial data, the normalizations $\tilde\varphi$ of the Chern-Ricci flow converge in $C^\infty$ to a solution of the equation $(\alpha+i\p\bar\p\varphi)^n=e^{\psi+c}\alpha^n$, for some constant $c$.
\end{corollary}

We note that there is a rich literature on Monge-Amp\`ere equations, including considerable progress using pluripotential theory. We refer to \cite{K, EGZ, DP, GZ,GZ2, PSS12, To1,To2, Na,Nb} and references therein.

\subsection{Hessian flows}

Hessian equations, where the Laplacian or the Monge-Amp\`ere determinant of the unknown function $u$ are replaced by the $k$-th symmetric polynomial of the eigenvalues of the Hessian of $u$, were introduced by Caffarelli, Nirenberg, and Spruck \cite{CNS}. More general right hand sides and K\"ahler versions were considered respectively by Chou and Wang \cite{CW} and Hou-Ma-Wu \cite{HMW}, who introduced in the process some of the key techniques for $C^2$ estimates that we discussed in \S 2. A general existence result on compact Hermitian manifolds was recently obtained by Dinew and Kolodziej \cite{DK},
Sun \cite{Su4}, and Sz\'ekelyhidi \cite{Sze}. See also Zhang \cite{Zh}. Again, we can derive this theorem as a corollary of Theorem \ref{theorem B}:

\begin{corollary}
\label{Cor4} 
Let $(X,\alpha)$ be a compact Hermitian $n$-dimensional manifold, and let $\chi$ be a positive real $(1,1)$-form which is $k$-positive for a given $k$, $1\leq k\leq n$. 
Consider the following parabolic flow for the unknown function $u$,
\bea
\label{kHessian parabolic}
\p_tu
=
\log{(\chi+i\p\bar\p u)^k\wedge \alpha^{n-k}\over\alpha^n}-\psi(z).
\eea
Then for any admissible initial data $u_0$, the flow admits a solution $u(z,t)$ for all time, and its normalization $\tilde u(z,t)$ converge in $C^\infty$
to a function $u_\infty\in C^\infty(X)$ so that $\o=\chi+i\p\bar\p u_\infty$ satisfies the following $k$-Hessian equation,
\bea
\label{kHessian}
\o^k\wedge \alpha^{n-k}=e^{\psi+c}\alpha^n.
\eea
\end{corollary}

\noindent
{\it Proof of Corollary \ref{Cor4}}.
This is an equation of the form (\ref{parabolic eq}), with $F=f(\lambda)=\log \,\sigma_k(\lambda)$, defined on the cone 
\bea
\Gamma_k=\{\lambda; \, \sigma_j(\lambda)>0,
\ j=1,\cdots,k\},
\eea
where ${{n}\choose{k}}\sigma_k$ is the $k$-th symmetric polynomial in the components $\lambda_j$, $1\leq j\leq n$. In our setting,
\bea
\sigma_k(\lambda[u])= {(\chi+i\p\bar\p u)^k\wedge \alpha^{n-k}\over\alpha^n}.
\eea
It follows from  \cite[Corollary 2.4]{Sp} that $g=\sigma_k^{1/ k}$ is concave and 
 $g_i=\frac{\p g}{\p \lambda_i}>0$  on $\Gamma_k$, hence $f=\log g$ satisfies the conditions (1-3) mentioned in \S 1.

The function $\underline u=0$ is a  subsolution of (\ref{kHessian parabolic})  and $f$ is in the unbounded case since for any $\mu=(\mu_1,\cdots,\mu_n)\in\Gamma_k$,
and any $1\leq i\leq n$,
\bea
{\rm lim}_{s\to\infty}\log\sigma_k(\mu_1,\cdots,\mu_i+s,\cdots,\mu_n)
=\infty.
\eea
 The desired statement follows then from Theorem \ref{theorem B}.

\subsection{The $J$ flow and quotient Hessian flows}

The $J$-flow on K\"ahler manifolds was introduced independently by Donaldson \cite{D1} and Chen \cite{C1}. The case $n=2$ was solved by Weinkove \cite{W1, W2}, and the case of general dimension by Song and Weinkove \cite{SW}, who identified a necessary and sufficient condition for the long-time existence and convergence of the flow as the existence of a K\"ahler form $\chi$ satisfying
\bea
\label{SW1}
nc\chi^{n-1}-(n-1)\chi^{n-2}\wedge \o>0
\eea
in the sense of positivity of $(n-1,n-1)$-forms. The constant $c$ is actually determined by cohomology. Their work was subsequently extended to inverse Hessian flows on K\"ahler manifolds by Fang, Lai, and Ma \cite{FLM}, and to inverse Hessian flows  on Hermitian manifolds  by Sun  \cite{Su1}.
These flows are all special cases of quotient Hessian flows on Hermitian manifolds. Their stationary points are given by the corresponding quotient Hessian equations. Our results can be applied to prove the following generalization to
quotient Hessian flows of the results of \cite{W1,W2,FLM}, as well as an alternative proof of a result of Sz\'ekelyhidi \cite[Proposition 22]{Sze} on the Hessian quotient equations. The flow (\ref{kl}) below has also been studied recently by Sun \cite{Su3} where he obtained a uniform $C^0$ estimate using Moser iteration. Our proof should be viewed as different from all of these, since its $C^0$ estimate uses neither Moser iteration nor strict $C^2$ estimates ${\rm Tr}_\alpha \chi_u\leq C\,e^{u-{\rm inf}_Xu}$.

\begin{corollary}
\label{Cor5}
Assume that $(X,\alpha)$ is a compact K\"ahler $n$-manifold, and fix $1\leq \ell<k\leq n$. Fix a closed $(1,1)$-form $\chi$ which is  $k$-positive, and assume that there exists a function $\underline u$ so that the form $\chi'=\chi+i\p\bar \p\underline{u}$ is closed $k$-positive and satisfies
\bea
\label{kl condition}
kc\,(\chi')^{k-1}\wedge \alpha^{n-k}
-\ell (\chi')^{\ell-1}\wedge \alpha^{n-\ell}>0
\eea
in the sense of the positivity of $(n-1,n-1)$-forms. Here $c=
{[\chi^\ell]\cup[\alpha^{n-\ell}]\over[\chi^k]\cup[\chi^{n-k}]}.$ Then for any admissible initial data $u_0\in C^\infty(X)$, the flow
\bea
\label{kl}
\p_t u=c-{\chi_u^\ell\wedge\alpha^{n-\ell}\over
\chi_u^k\wedge \alpha^{n-k}}
\eea
admits a solution $u$ for all time, and it converges to a smooth function $u_\infty$. The
form $\o=\chi+i\p\bar\p u_\infty$ is $k$-positive and satisfies the equation
\bea
\label{quotient}
\o^\ell\wedge \alpha^{n-\ell}=c\,\o^k\wedge\alpha^{n-k}.
\eea
\end{corollary}

\medskip
\noindent
{\it Proof of Corollary \ref{Cor5}}. 
The flow (\ref{kl}) is of the form (\ref{parabolic eq}), with $$ f(\lambda)=- \frac{\sigma_\ell (\lambda)}{ \sigma_k(\lambda)},$$
 defined on the cone 
\bea
\Gamma_k=\{\lambda; \, \sigma_j(\lambda)>0,
\ j=1,\cdots,k\}.
\eea
By the Maclaurin's inequality (cf. \cite{Sp}), we have  $\sigma_k^{1/k}\leq \sigma_\ell^ {1/\ell}$ on $\Gamma_k$, hence $f(\lambda)\rightarrow -\infty$ as $\lambda\rightarrow \p \Gamma_k$.
It follows from \cite[Theorem 2.16]{Sp} that the function $g=(\sigma_k/\sigma_\ell)^{1\over(k-\ell)}$ satisfies $g_i=\frac{\p g}{\p \lambda_i}>0,$ $\forall i=1,\ldots,n$ and $g$ is concave on $\Gamma_k$. Therefore $f=-g^{-(k-\ell)}$ satisfies the condition (1-3) spelled out in \S 1. Moreover, $f$ is in the bounded case with
$$
f_\infty(\lambda')= - { \ell\sigma_{\ell-1}(\lambda')\over k\sigma_{k-1}(\lambda' }\quad  {\rm where} \quad \lambda'\in \Gamma_\infty=\Gamma_{k-1}.
$$
 We can assume that $u_0=0$ by replacing $\chi$ (resp. $u$ and $\underline{u}$) by $\chi+i\p\bar \p u_0$ (resp. $u-u_0$ and $\underline{u}-u_0$). 
The inequality (\ref{kl condition})  infers that $\underline{u}$ is a subsolution of the equation (\ref{kl}). Indeed,
for any $(z,t)\in X\times[0,\infty)$,  set $\mu=\lambda(B)$, $B^i{}_j=\alpha^{j\bar k}(\chi_{\bar kj}+ \underline u_{\bar k j})(z,t)$. Since  $\underline{u}$ is independent of $t$, it follows from Lemma \ref{criterion} and the symmetry of $f$ that we just need to show that  for any $z\in X$  if
$\mu'=(\mu_1,\cdots,\mu_{n-1})$ then
\bea
 {\rm lim}_{s\to\infty}f(\mu' , \mu_n +s)>-c.
\eea
This means
\bea
\label{ineq_subsolution}
f_\infty(\mu' )= - { \ell\sigma_{\ell-1}(\mu')\over k\sigma_{k-1}(\mu' )} >-c.
\eea
As in \cite{Sze}, we restrict to the tangent space of $X$ spanned by by the eigenvalues corresponding to $\mu'$. Then on this subspace
\bea
\sigma_j(\mu')={\chi^{j-1}\wedge \alpha^{n-j}\over \alpha^{n-1}}
\eea
for all $j$. Thus the preceding inequality  is equivalent to
\bea
kc(\chi')^k\wedge \alpha^{n-k}-\ell (\chi')^{\ell-1}\wedge \alpha^{n-\ell}>0.
\eea
By a priori estimates in Section 2, the solution exists for all times. We now use the second statement in  Lemma \ref{convergence lemma} to prove the convergence.  It suffices to check that $u$ is uniformly bounded  in $X\times [0,+\infty)$ and for all $t>0$, there exists $y$ such that $\p_t u(y,t)=0$. The second condition is straightforward since 
$$
\int_X \p_t u \chi_u^{k}\wedge \alpha^{n-k} =0.
$$
For the uniform bound we make use of the following lemma
\begin{lemma}
\label{I_functional}
Let  $\phi\in C^\infty(X)$ function and  $\{\varphi_s\}_{s\in [0,1]}$ be a path with  $\varphi(0)= 0$ and $ \varphi(1)=\phi$.  Then we have
\bea
\label{I_k}
\int^1_0\int_X  \frac{\partial \varphi }{\partial s}\chi_{\varphi }^k\wedge\alpha^{n-k} ds=\frac{1}{k+1}\sum_{j=0}^k \int_X \f \chi_\f^j\wedge\chi^{k-j} \wedge\alpha^{n-k},
\eea
so the left hand side is independent of $\varphi$. 
Therefore we can define the following functional
\begin{equation}
\label{functional}
I_k(\phi)=\int^1_0\int_X  \frac{\partial \varphi }{\partial s}\chi_{\varphi }^k\wedge\alpha^{n-k} ds.
\end{equation}
\end{lemma}
We remark that when $k=n$ and $\chi$ is K\"ahler, this functional is  well-known (see for instance \cite{W2}). We discuss here the general case.

\medskip
\noindent
{\it Proof of Lemma \ref{I_functional}.} 
Observe that 
  \begin{eqnarray}
\label{func_0}
\int^1_0\int_X  \frac{\partial \varphi }{\partial s}\chi_{\varphi }^k\wedge\alpha^{n-k} ds =\sum_{j=1}^{k} {{k}\choose{j}}\int^1_0\int_X  \frac{\partial \varphi }{\partial s}(i\p\bar{\p}  \varphi)^j \wedge \chi^{k-j}\wedge\alpha^{n-k} ds. 
\end{eqnarray} For any $j=0,\ldots,k $ we have
\begin{eqnarray} \nonumber
\int^1_0\int_X  \frac{\partial \varphi }{\partial s}(i\p\bar{\p} \varphi)^j \wedge \chi^{k-j}\wedge\alpha^{n-k} ds 
&=&  \int_0^1 \frac{d}{ds} \left(\int_X \varphi  (i\p\bar{\p}  \varphi)^j \wedge \chi^{k-j}\wedge\alpha^{n-k} \right)ds\\ \nonumber
&& -  \int_0^1 \int_X  \varphi \frac{\p}{\p s}\left( (i\p\bar{\p}  \varphi)^j \wedge \chi^{k-j}\wedge\alpha^{n-k} \right)ds\\ \label{func_1}
&=&\int_X \phi  (i\p\bar{\p}  \phi)^j \wedge \chi^{k-j}\wedge\alpha^{n-k} \\\nonumber
&& -  \int_0^1 \int_X  \varphi \frac{\p}{\p s}\left( (i\p\bar{\p}  \varphi)^j \wedge \chi^{k-j}\wedge\alpha^{n-k} \right)ds
\end{eqnarray} 
We also have
\begin{eqnarray}\nonumber
 \int_0^1 \int_X  \varphi \frac{\p}{\p s}\left( (i\p\bar{\p}  \varphi)^j \wedge \chi^{k-j}\wedge\alpha^{n-k} \right)ds&=&\int_0^1 \int_X j  \varphi \left( i\p\bar{\p} \frac{\p\varphi}{\p s}\right) \wedge(i\p\bar{\p}  \varphi)^{j-1} \wedge \chi^{k-j}\wedge\alpha^{n-k}ds\\\label{func_2}
 &=&\int_0^1 \int_X j   \frac{\partial \varphi }{\partial s}(i\p\bar{\p}  \varphi)^j \wedge \chi^{k-j}\wedge\alpha^{n-k} ds,
\end{eqnarray}
here we used in the second identity the integration by parts and the fact that $\chi$ and $\alpha$  are closed.
Combining (\ref{func_1}) and (\ref{func_2})  yields 
\begin{equation}
\int^1_0\int_X  \frac{\partial \varphi }{\partial s}(i\p\bar{\p}  \varphi)^j \wedge \chi^{k-j}\wedge\alpha^{n-k} ds=\frac{1}{j+1} \int_X \phi  (i\p\bar{\p}  \phi)^j \wedge \chi^{k-j}\wedge\alpha^{n-k}.
\end{equation}
Therefore (\ref{func_0}) implies that
\begin{eqnarray} \nonumber
\int^1_0\int_X  \frac{\partial \varphi }{\partial s}\chi_{\varphi }^k\wedge\alpha^{n-k} ds &=&\sum_{j=1}^{k} {{k}\choose{j}}\frac{1}{j+1} \int_X \phi   (i\p\bar{\p} \phi)^j \wedge \chi^{k-j}\wedge\alpha^{n-k}\\\nonumber
&=&\sum_{j=1}^{k} {{k}\choose{j}}\frac{1}{j+1} \int_X \phi   (\chi_\phi-\chi)^j \wedge \chi^{k-j}\wedge\alpha^{n-k}\\ \label{func_4}
&=& \sum_{j=1}^{k} {{k}\choose{j}}\frac{1}{j+1} \int_X  \sum_{p=0}^{j}{{j}\choose{p}} (-1)^{j-p} \phi \chi_\phi^p \wedge\chi^{k-p}\wedge \alpha^{n-k}\\ \nonumber
&=& \sum_{p=0}^{k}\left( \sum_{j=p}^k  {{k}\choose{j}}\frac{1}{j+1}{{j}\choose{p}} (-1)^{j-p}\right)  \int_X\phi \chi_\phi^p \wedge\chi^{k-p}\wedge \alpha^{n-k}.
 \end{eqnarray} 
By changing $m=j-p$, we get
 \begin{eqnarray}
\label{func_5}
 \sum_{j=p}^k  {{k}\choose{j}}\frac{1}{j+1}{{j}\choose{p}} (-1)^{j-p}&=&{{k}\choose{p}}  \sum_{m=0}^{k-p} \frac{(-1)^m}{m+p+1} {{k-p}\choose{m}}.
 \end{eqnarray}
The right hand side can be computed by 
\begin{eqnarray*}
{{k}\choose{p}}  \sum_{m=0}^{k-p} \frac{(-1)^m}{m+p+1} {{k-p}\choose{m}} &=&{{k}\choose{p}} \int_0^1(1-x)^{k-p}x^pdx\\
&=&{{k}\choose{p}} p!\int_{0}^1 \frac{1}{(k-p+1)\ldots k} (1-x)^kdx=
\frac{1}{k+1},
\end{eqnarray*}
where we used the integration by parts $p$ times in the second identity. Combining this with (\ref{func_4}) and (\ref{func_5}) we get the desired identity (\ref{I_k}). Q.E.D.

\medskip
 We now have for any $t^*>0$, along the flow
\begin{eqnarray*}
I_k(u(t^*))=\int_0^{t^*}\int_X \frac{\partial u}{\partial t}\chi_u^k\wedge\alpha^{n-k}
=\int_0^{t^*}  \left(c- \frac{\chi_u^\ell\wedge\alpha^{n-\ell}}{\chi_u^k	\wedge\alpha^{n-k}}\right)\chi_u^k\wedge\alpha^{n-k}=0.
\end{eqnarray*}
As in Weinkove \cite{W1,W2}, 
there exist $C_1,C_2>0$ such that for all $t\in [0,\infty)$, 
\begin{equation}\label{harnack_ineq_kl}
0\leq \sup_X u(.,t)\leq -C_1 \inf_X u(.,t)+C_2.
\end{equation}
Indeed, in view of (\ref{I_k}),  $I_k(u)=0$ along the flow implies that 
\begin{equation}\label{identity_u}
\sum_{j=0}^{k} \int_X u\chi_{u}^j\wedge \chi^{k-j}\wedge\alpha^{n-k}=0,
\end{equation}
hence
 $\sup_X u\geq 0$ and $\inf_X u\leq 0$. 
\medskip
For the right inequality in (\ref{harnack_ineq_kl}), we remark that there exists a positive constant $B$ such that
$$\alpha^n\leq B\chi^k\wedge \alpha^{n-k}.$$
Therefore combining with (\ref{identity_u}) gives
\begin{eqnarray*}
\int_X u\alpha^n &=& \int_X(u-\inf_X  u)\alpha^n+\int_X \inf_Xu \,\alpha^n\\
&\leq& B \int_X(u-\inf_X u)\chi^k \wedge\alpha^{n-k} +\inf_X u\int_X \alpha^n\\
&= & -B \sum_{j=1}^{k}\int_Xu\chi_{u}^j\wedge \chi^{k-j}\wedge \alpha^{n-k} +\inf_Xu\left(\int_X \alpha^n-B\int_X \chi^{k}\wedge \alpha^{n-k} \right) \\
&= & -B \sum_{j=1}^{k}\int_X\left( u-\inf_X u\right) \chi_{u}^j\wedge \chi^{k-j}\wedge \alpha^{n-k} +\inf_Xu\left(\int_X \alpha^n-B(k+1)\int_X \chi^{k}\wedge \alpha^{n-k} \right)\\
&\leq & \inf_X u\left(\int_X \alpha^n-B(k+1)\int_X \chi^k\wedge \alpha^{n-k} \right)=-C_1\inf_X u. 
\end{eqnarray*}
Since $\Delta_\alpha u\geq -\tr_\alpha \chi\geq -A$, using the fact that the  Green's function  $G(.,.)$ of $\alpha$ is   bounded from below we infer that 
\begin{eqnarray*}
u(x,t)&= &\int_X u\alpha^n-\int_X \Delta_\alpha  u(y,t)G(x,y)\alpha^n(y)\\
&\leq&  - C_1\inf_X u + C_2.
\end{eqnarray*}
Hence we obtain the Harnack inequality,  $\sup_X  u \leq-C_1\inf_X u+C_2  $.

\medskip
Since we can normalize $\underline{u}$ by $\sup_X \underline{u}=0$, the left inequality in (\ref{identity_u}) implies
$$\sup_X(u(\cdot, t)-\underline{u}(\cdot, t))\geq 0.$$
It follows from Lemma \ref{C0} that $$u\geq \underline{u}-C_3$$
for some constant $C_3$. This give a lower bound for $u$ since $\underline{u}$ is bounded.  The Harnack inequality in (\ref{harnack_ineq_kl}) implies then a uniform bound for $u$. Now the second statement in Lemma \ref{convergence lemma} implies the convergence of $u$. Q.E.D.

\bigskip
A natural generalization of the Hessian quotient flows on Hermitian manifolds is the following flow
\bea
\label{log}
\p_tu=\log {\chi_u^k\wedge \alpha^{n-k}\over\chi_u^\ell\wedge \alpha^{n-\ell}}-
\psi
\eea
where $\psi\in C^\infty(X)$, the admissible cone is $\Gamma_k$, $1\leq \ell<k\leq n$, and $\chi_u=\chi+i\p\bar\p u$. This flow was introduced by Sun \cite{Su1} when $k=n$. We can apply Theorem \ref{theorem C} to obtain the following result, which is analogous to one of the main results in Sun \cite{Su1}, and analogous to the results of Song-Weinkove \cite{SW} and Fang-Lai-Ma \cite{FLM} for $k=n$:

\begin{corollary}
\label{Cor6}
Let $(X,\alpha)$ be a compact Hermitian manifold and $\chi$ be a $(1,1)$-form  which is  $k$-positive.
Assume that there exists a form $\chi'=\chi+i\p\bar\p \underline{u}$ which is $k$-positive, and satisfies
\bea
k\,(\chi')^{k-1}\wedge \alpha^{n-k}
-e^\psi\,\ell (\chi')^{\ell-1}\wedge \alpha^{n-\ell}>0
\eea
in the sense of the positivity of $(n-1,n-1)$-forms. Assume further that there exists an admissible $u_0\in C^\infty(X)$ satisfying
\bea
\label{ineq_condition}
e^\psi\geq {\chi^k_{u_0}\wedge\alpha^{n-k}\over \chi_{u_0}^\ell\wedge \alpha^{n-\ell}}
\eea
Then the flow (\ref{log}) 
admits a smooth solution for all time with initial data $u_0$. Furthermore,
there exists a unique constant $c$ so that the normalization
\bea
\tilde u=u-{1\over [\alpha^n]}\int_X u\alpha^n
\eea
converges in $C^\infty$ to a function $u_\infty$ with
$\o_\infty=\chi+i\p\bar\p u_\infty$ satisfying
\bea
\o_\infty^k \wedge \alpha^{n-k}=e^{\psi+c}\o_\infty^\ell\wedge\alpha^{n-\ell}.
\eea
\end{corollary}

\noindent
{\it Proof of Corollary \ref{Cor6}}. This equation is of the form (\ref{parabolic eq}), with
\bea
F(A)=f(\lambda)=\log {\sigma_k(\lambda)\over\sigma_\ell(\lambda)}, \ 
{\rm with}\ \lambda=\lambda(A),
\eea
defined on $\Gamma_k$.
As in the proof of Corollary \ref{Cor5}  we also have that  $f$ satisfies the conditions (1-3) mentioned in \S 1.  Moreover, $f$ is in the bounded case with
$$
f_\infty(\lambda')=\log{ k\sigma_{k-1}(\lambda')\over\ell\sigma_{\ell-1}(\lambda')} \quad {\rm where} \quad \lambda'\in \Gamma_\infty=\Gamma_{k-1}.
$$
It suffices to verify that $\underline u=0$ is a subsolution of the equation (\ref{log}).
For any $(z,t)\in X\times[0,\infty)$,  set $\mu=\lambda(B)$, $B^i{}_j=\alpha^{j\bar k}\chi_{\bar kj}(z,t)$. Since  $\underline{u}$ is independent of $t$, Lemma \ref{criterion} implies that we just need to show that  for any $z\in X$   if
$\mu'=(\mu_1 ,\cdots,\mu_{n-1})$,
\bea
{\rm lim}_{s\to\infty}f(\mu' ,s)>\psi(z).
\eea
This means 
\bea
f_\infty(\mu') =\log{ k\sigma_{k-1}(\mu')\over\ell\sigma_{\ell-1}(\mu')}>\psi(z),
\eea
where we restrict to the tangent space of $X$ spanned by by the eigenvalues corresponding to $\mu'$. As the argument in the proof of Corollary \ref{Cor5}, this inequality is  equivalent to
\bea
k\chi^k\wedge \alpha^{n-k}-\ell e^\psi \chi^{\ell-1}\wedge \alpha^{n-\ell}>0.
\eea
Moreover, the condition (\ref{ineq_condition})  is equivalent to
\bea
0= \underline{u}\geq F(A[u_0])-\psi. 
\eea
We can now apply Theorem \ref{theorem C} to complete the proof. Q.E.D

\medskip
In the case of $(X,\alpha)$ compact K\"ahler, the condition on $\psi$ can be simplified, and we obtain an alternative proof to the main result of Sun in \cite{Su2}.
We recently learnt that Sun \cite{Su6} also provided independently another proof of 
\cite{Su2} using the same flow as below:

\begin{corollary}
\label{Cor7}
Let $(X,\alpha)$ be K\"ahler and $\chi$ be a  $k$-positive closed $(1,1)$-form. Assume that there exists a closed form $\chi'=\chi+i\p\bar\p \underline{u}$ which is $k$-positive, and satisfies
\bea
\label{log condition 2}
k\,(\chi')^{k-1}\wedge \alpha^{n-k}
-e^\psi\,\ell (\chi')^{\ell-1}\wedge \alpha^{n-\ell}>0
\eea
in the sense of the positivity of $(n-1,n-1)$-forms. Assume further that 
\bea
\label{log  condition 3}
e^\psi\geq c_{k,\ell}=
{[\chi^k]\cup[\chi^{n-k}]\over [\chi^\ell]\cup[\alpha^{n-\ell}]}.
\eea
Then for any admissible initial data $u_0\in C^\infty(X)$, the flow (\ref{log}) 
admits a smooth solution for all time. Furthermore,
there exists a unique constant $c$ so that the normalization
\bea
\tilde u=u-{1\over [\alpha^n]}\int_X u\alpha^n
\eea
converges in $C^\infty$ to a function $u_\infty$ with
$\o_\infty=\chi+i\p\bar\p u_\infty$ satisfying
\bea
\o_\infty^k \wedge \alpha^{n-k}=e^{\psi+c}\o_\infty^\ell\wedge\alpha^{n-\ell}.
\eea
\end{corollary}

\noindent
{\it Proof of Corollary \ref{Cor6}}. By the same argument above, the admissible function $\underline{u}\in C^\infty(X)$  with $\sup_X \underline{u}=0$ satisfying (\ref{log condition 2})
 is a $C$-subsolution. As explained in the proof of Corollary \ref{Cor5}, we can assume  that $u_0=0$.
 
\medskip 
We first observe that along the flow, the functional  $I_\ell$ defined in Lemma \ref{I_functional} is decreasing. Indeed, using Jensen's inequality and then (\ref{log  condition 3}) we have 
\begin{eqnarray}
\frac{d}{dt}I_\ell (u)&=&\int_X \frac{\partial u}{\partial t} \chi_u^{\ell}\wedge \alpha^{n-\ell} 
=\int_X \left(\log \frac{ \chi_u^k\wedge \alpha^{n-k}}{\chi_u^\ell\wedge \alpha^{n-\ell}} -\psi \right)\chi_u^{\ell}\wedge \alpha^{n-\ell}\nonumber\\
&\leq& \log c_{k,\ell}\int_X \chi_u^\ell\wedge \alpha^{n-\ell}  -\int_X\psi \chi_u^\ell\wedge \alpha^{n-\ell} 
\leq 0.
\end{eqnarray}
Set
\bea
\hat{u}:= u-h(t),
\qquad
h(t)=\frac{I_\ell(u)}{\int_X \chi^\ell\wedge \alpha^{n-\ell}}.
\eea
For any $ t^* \in [0,\infty)$ we have
\begin{eqnarray}
I_\ell(\hat u( t^*))=\int_0^{ t^*}\int_X \frac{\partial \hat u}{\partial  t} \chi_u^\ell\wedge\alpha^{n-\ell}
=\int_0^{t^*} \int_X  \left( \frac{\partial u }{\partial t} -\frac{ 1}{\int_X \chi^\ell\wedge \alpha^{n-\ell}}\frac{d}{dt}I_\ell(u) \right)\chi_u^\ell\wedge \alpha^{n-\ell}=0.
\nonumber
\end{eqnarray}
By the same argument in Corollary \ref{Cor5}, we deduce that there exist $C_1,C_2>0$ such that 
\begin{equation}\label{bound_suf_inf_3}
0\leq \sup_X\hat u(.,t)\leq -C_1 \inf_X \hat{u}(.,t)+C_2,
\end{equation}
for all $t\in [0,\infty)$. By our choice, $\sup_X \underline u=0$, and (\ref{bound_suf_inf_3}) implies that $$\sup_X(u-h(t)-\underline{u} )=\sup_X( \hat{u}-\underline u) \geq 0,\, \forall t \geq 0.$$ 
Since  $ I_\ell(u)$ is  decreasing along the flow, we also have $h'(t)
\leq 0$.   Theorem \ref{theorem C} now gives us the required result. Q.E.D.

\medskip
Similarly, we can consider the flow (\ref{parabolic eq}) with
\bea
\label{inverse Hessian flow}
\p_tu
=-\left({\chi_u^\ell \wedge \alpha^{n-\ell}\over \chi_u^k\wedge \alpha^{n-k}}\right)^{1\over k-\ell}+\psi(z), \quad u(z,0)=0,
\eea
where $1\leq \ell<k\leq n$.
When $(X,\alpha)$ is K\"ahler, $\psi$ is constant and  $k=n$, this is the inverse Hessian flow studied by Fang-Lai-Ma \cite{FLM}. We can apply Theorem \ref{theorem C} to obtain another corollary which  is  analogous to the main result of Fang-Lai-Ma \cite{FLM}.

\begin{corollary}
\label{Cor8}
Let $(X,\alpha)$, and $\chi$ as in Corollary \ref{Cor6}. Assume further that $\psi \in C^\infty(X,\R^+)$ and there exists a smooth function $\underline u$ with $\chi'=\chi+i\p\bar\p \underline{u}$ a $k$-positive $(1,1)$-form which satisfies
\bea
\label{condition_8}
k \psi^{k-\ell} (\chi')^{k-1}\wedge \alpha^{n-k}-\ell (\chi')^{n-\ell-1}\wedge\alpha^{n-\ell}>0
\eea
in the sense of positivity of $(n-1,n-1)$ forms, and
\bea
\label{inverse_condition1}
\psi^{k-\ell}\leq \frac{\chi^{\ell}\wedge \alpha^{n-\ell}}{\chi^k\wedge \alpha^{n-k}}.
\eea
Then the flow (\ref{inverse Hessian flow}) exists for all time, and there is a unique constant $c$ so that the normalized function $\tilde u$ converges to a function $u_\infty$ with
$\omega=\chi+i\p\bar\p u_\infty$ a $k$-positive form satisfying the equation
\bea
\omega^{n-\ell}\wedge \alpha^{n-\ell}
=
(\psi+c)^{k-\ell} \o^{k}\wedge\alpha^{n-k}.
\eea
In particular, if $(X,\alpha)$ is K\"ahler, we assume further that $\chi$ is closed,  then the condition (\ref{inverse_condition1}) can be simplified as
\begin{equation}
\label{inverse_condition2}
\psi^{k-\ell}\leq c_{\ell,k}=\frac{[\chi^{\ell}]\cup [\alpha^{n-\ell}]}{[\chi^k]\cup [\alpha^{n-k}]}.
\end{equation}
\end{corollary}

\noindent
{\it Proof of Corollary \ref{Cor8}}. This equation is of the form (\ref{parabolic eq}), with
\bea
F(A)=f(\lambda)=- \left(\frac{\sigma_\ell(\lambda)}{\sigma_k(\lambda)}\right)^{1\over k-\ell}, \ 
{\rm with}\ \lambda=\lambda(A),
\eea
defined on $\Gamma_k$. 
As in Corollary \ref{Cor5}, it follows from the Maclaurin's inequality, the monotonicity and concavity of $g=(\sigma_k/\sigma_\ell)^{1\over k-\ell}$  (cf. \cite{Sp}) that  $f$ satisfies the conditions (1-3) spelled out in \S 1.
Moreover, $f$ is in the bounded case with
$$
f_\infty(\lambda')=- \left(\frac{\ell\sigma_{\ell-1}(\lambda')}{k\sigma_{k-1}(\lambda')}\right)^{1\over k-\ell} \quad {\rm where}\quad \lambda'\in \Gamma_\infty=\Gamma_{k-1}. 
$$
 In addition, as the same argument in previous corollaries, the condition (\ref{condition_8})  is equivalent to that $\underline{u}=0$ is a $C$-subsolution for (\ref{inverse Hessian flow}). Moreover, the condition (\ref{inverse_condition1})  implies that 
\bea
0= \underline{u}\geq F(A[0])+\psi. 
\eea
We can now apply Theorem \ref{theorem C} to get the first result.

\medskip
Next, assume that $(X,\alpha)$ is K\"ahler  and $\chi$ is closed. As in Corollary \ref{Cor7} and \cite{FLM}, the functional $I_\ell$ (see Lemma \ref{I_functional}) is decreasing along the flow.  Indeed, using (\ref{inverse_condition2}),
\begin{eqnarray}
\label{d/dt I}
\frac{d}{dt}I_\ell (u)&=&\int_X \frac{\partial u}{\partial t} \chi_u^{\ell}\wedge \alpha^{n-\ell} 
=\int_X \left(- \left(\frac{\sigma_\ell(\lambda)}{\sigma_k(\lambda)}\right)^{1\over k-\ell}+\psi \right)\chi_u^{\ell}\wedge \alpha^{n-\ell}\nonumber\\
&\leq & -\int_X \left(\frac{\sigma_\ell(\lambda)}{\sigma_k(\lambda)}\right)^{1\over k-\ell} \chi_u^{\ell}\wedge \alpha^{n-\ell} +c_{\ell,k}^{1\over k-\ell}\int_X \chi_u^{\ell}\wedge \alpha^{n-\ell}.
\end{eqnarray}
Using  the H\"older inequality, we get
\begin{eqnarray*}
\int_X \chi_u^{\ell}\wedge \alpha^{n-\ell} &= &\int_X \sigma_\ell \,\alpha^n=\int_X\left( {\sigma_\ell\over \sigma_k^{1/( k-\ell+1)} }\right) \sigma_k^{1\over k-\ell+1}\,\alpha^n\\
&\leq  & \left[  \int_X\left( {\sigma_\ell\over \sigma_k^{1/( k-\ell+1)} }\right)^{k-\ell+1\over k-\ell}  \alpha^n\right]^{k-\ell\over k-\ell +1}  \left(\int_ X \sigma_k\,\alpha^n\right)^{1\over k-\ell+1}\\
&=&\left[ \int_X \left(\frac{\sigma_\ell(\lambda)}{\sigma_k(\lambda)}\right)^{1\over k-\ell} \chi_u^{\ell}\wedge \alpha^{n-\ell} \right]^{k-\ell\over k-\ell +1}   \left(\int_ X \chi_u^{k}\wedge \alpha^{n-k}\right)^{1\over k-\ell+1}\\
&= &\left[ \int_X \left(\frac{\sigma_\ell(\lambda)}{\sigma_k(\lambda)}\right)^{1\over k-\ell} \chi_u^{\ell}\wedge \alpha^{n-\ell} \right]^{k-\ell\over k-\ell +1}   c_{\ell,k}^{-1\over k-\ell +1}\left(\int_ X \chi_u^{\ell}\wedge \alpha^{n-\ell}\right)^{1\over k-\ell+1}.
\end{eqnarray*}
This implies that $$c_{\ell,k}^{1\over k-\ell}\int_X \chi_u^{\ell}\wedge \alpha^{n-\ell}\leq \int_X \left(\frac{\sigma_\ell(\lambda)}{\sigma_k(\lambda)}\right)^{1\over k-\ell} \chi_u^{\ell}\wedge \alpha^{n-\ell}, $$
hence
$d I_\ell (u)/ dt\leq 0$.

\medskip
 For the rest of the proof, we follow the argument in Corollary  \ref{Cor7}, starting from the fact that  $I_\ell (\hat u)=0$ where
 $$\hat u=u-\frac{I_\ell (u)}{\int_X \chi^\ell \wedge \alpha^{n-\ell}}.$$
Then we obtain the Harnack inequality  \begin{equation}
0\leq \sup_X\hat u(.,t)\leq -C_1 \inf_X \hat{u}(.,t)+C_2,
\end{equation} 
for some constants $C_1,C_2>0$.  Finally, Theorem \ref{theorem C}  gives us the last claim. Q.E.D.

\subsection{Flows with mixed Hessians $\sigma_k$}
Our method can be applied to solve other equations containing many terms of $\sigma_k$.  We illustrate this with  the  equation
\begin{equation}\label{general_sigma_k}
\sum_{j=1}^\ell c_j \chi_u^j\wedge \alpha^{n-j}=c\chi_u^k\wedge \alpha^{n-k}
\end{equation}
on a K\"ahler manifold $(X,\alpha)$,
where $1\leq \ell<k\leq n$, $c_j\geq 0$ are given non-negative constants, and $c\geq 0$ is determined by $c_j$  by integrating the equation over $X$. 

\medskip
When $k=n$, It was conjectured by  Fang-Lai-Ma \cite{FLM} that this equation is solvable assuming that 
$$
nc\chi'^{n-1}-\sum_{k=1}^{n-1}kc_k\chi'^{k-1}\wedge \alpha^{n-k}>0,
$$
for some closed $k$-positive form $\chi'=\chi+i\p\bar\p v$. This conjecture was solved recently by Collins-Sz\'ekelyhidi \cite{CSze} using the continuity method. 
An alternative proof by flow methods is in Sun \cite{Su5}.
Theorem \ref{Cor9} stated earlier in the Introduction is an existence result 
for more general equations (\ref{general_sigma_k}) using  the flow (\ref{combination})
In particular, it gives a parabolic proof of a generalization of the conjecture due to Fang-Lai-Ma \cite[Conjecture 5.1]{FLM}.  We also remark that the flow (\ref{combination}) was mentioned in Sun \cite{Su1}, but no result given there, to the best of our understanding.

\medskip

{\it Proof of Theorem \ref{Cor9}}.  This equation is of the form (\ref{parabolic eq}), with 
$$F(A)=f(\lambda )= -{\sum_{j=1}^\ell c_j\sigma_j(\lambda)\over \sigma_k(\lambda)}+c,$$
defined on the cone $\Gamma_k$. As in the proof of Corollary \ref{Cor5}, for any $j=1,\ldots,\ell$, the function $-\sigma_j/\sigma_k$ on $\Gamma_k$  satisfies the conditions (1-3) in \S 1, so does $f$.  We also have that $f$ is in the bounded case with
$$
f_\infty(\lambda')=-{\sum_{j=1}^\ell jc_j\sigma_{j-1}(\lambda)\over k\sigma_{k-1}(\lambda)} \quad {\rm where}\quad \lambda'\in \Gamma_\infty=\Gamma_k.
$$

\medskip
\noindent
Suppose $\chi'=\chi+i\p\bar \p \underline{u}$ with $\sup_X \underline{u}=0$ satisfies 
$$kc(\chi')^{k-1}\wedge\alpha^{n-k}-\sum_{j=1}^{\ell}jc_j(\chi')^{j-1}\wedge \alpha^{n-j}>0.$$
By the same argument in Corollary \ref{Cor5}, this is equivalent to that  $\underline u$ is a $C$-subsolution of  (\ref{combination}).  Observe that for all $t^{*}>0$,
\begin{eqnarray}
I_k(u (t^*))&=&\int^{t^*}_0 \int_X \frac{\partial u}{\partial t}\chi_u^k\wedge\alpha^{n-k}=\int^{t^*}_0 \int_X \left( c-\frac{\sum_{j=1}^\ell c_j \sigma_j(\lambda) }{\sigma_k (\lambda)} \right)\chi_u^k\wedge\alpha^{n-k}\nonumber\\
&=& \int^{t^*}_{0}\left(c\int_X \chi_u^k\wedge\alpha^{n-k} -\sum_{j=1}^\ell c_j \int_X\chi_u^j\wedge \alpha^{n-j}  \right) 
=0.
\end{eqnarray}
Therefore Lemma \ref{I_functional} implies that
$$\sum_{j=0}^{k} \int_X u\chi_u^j\wedge \chi^{k-j}\wedge\alpha^{n-k}=0.$$
Therefore we can obtain the Harnack inequality as in Corollary  \ref{Cor5}:
\bea
0\leq \sup_X u(.,t)\leq -C_1 \inf u(.,t)+C_2,
\eea
and $\inf_X u<0$, for some positive constants $C_1,C_2$. Lemma \ref{C0} then gives a uniform bound for $u$. Since $$\int_X \partial_t u\chi_u^k\wedge\alpha^{n-k}=0,$$
 for any $t>0$, there exists $y=y(t)$ such that $\partial_t u(y,t)=0$. The rest of the proof is the same to the proof of Corollary \ref{Cor5} where we used Lemma \ref{convergence lemma} to imply the convergence of the flow. Q.E.D.
 
 \medskip
We observe that equations mixing several Hessians seem to appear increasingly frequently in complex geometry. A recent example of particular interest is the Fu-Yau equation \cite{FY1,FY2, PPZ1, PPZ4} and its corresponding geometric flows \cite{PPZ5}.

\subsection{Concluding Remarks}

We conclude with a few open questions. 

\smallskip

It has been conjectured by Lejmi and Sz\'ekelyhidi \cite{LSze} that conditions of the form (\ref{SW1}) and their generalizations can be interpreted as geometric stability conditions. This conjecture has been proved in the case of the $J$-flow on toric varieties  by Collins and Sz\'ekelyhidi \cite{CSze}. Presumably there should be similar interpretations in terms of stability of the conditions formulated in the previous section. A discussion of stability conditions for constant scalar curvature K\"ahler metrics can be found in \cite{PS10}.

\smallskip
It would also be very helpful to have a suitable geometric interpretation of conditions such as the one on the initial data $u_0$. Geometric flows whose behavior may behave very differently depending on the initial data include the anomaly flows studied in \cite{PPZ2}, \cite{PPZ7}, \cite{FHP}.

\smallskip
For many geometric applications, it would be desirable to extend the theory of subsolutions to allow the forms $\chi$ and $\psi$ to depend on time as well as on $u$ and $\nabla u$.

\medskip
\noindent
{\bf Acknowledgements}: The second-named author would like to thank his thesis advisor, Vincent Guedj for his constant support and encouragement, and also Yuxin Ge for some useful discussions. This work was begun when the second-named author was visiting the department of Mathematics of Columbia University, and completed when the first-named author was visiting the Institut de Mathematiques of Toulouse. Both authors would like to thank these institutions for their hospitality.
They would also like to thank Nguyen Van Hoang for a careful reading of the paper.

{\it Email addresses:}

phong@math.columbia.edu,
Tat-Dat.To@math.univ-toulouse.fr

\end{document}